\documentclass{amsart}

\usepackage[latin1]{inputenc}

\usepackage{amsmath,amssymb,amscd,latexsym}

\usepackage[dvips]{graphicx}

\usepackage[english]{babel}

\newcommand{\bor}[1][]{\mathcal{B}(#1)}

\newcommand{\communique}{\leftrightarrow}
\newcommand{\N}{\mathbb{Z}_{+}}

\newcommand{\Z}{\mathbb{Z}}

\newcommand{\Zd}{\mathbb{Z}^d}

\newcommand{\Q}{\mathbb{Q}}
\newcommand{\Qjoli}{\mathcal{Q}}
\newcommand{\Pcond}{\overline{\mathbb{P}}}
\newcommand{\R}{\mathbb{R}}
\newcommand{\Rd}{\mathbb{R}^d}

\renewcommand{\P}{\mathbb{P}}
\newcommand{\E}{\mathbb{E}\ }
\newcommand{\F}{\E_{\Pcond}\ }

\newcommand{\Ed}{\mathbb{E}^d}
\newcommand{\Edeux}{\mathbb{E}^2}

\newcommand{\Ber}{\text{Ber}}

\newcommand{\as}{\text{ a.s.}}

\renewcommand{\epsilon}{\varepsilon}
\renewcommand{\limsup}{\overline{\lim}}
\renewcommand{\liminf}{\underline{\lim}}

\newcommand{\ie}{\emph{i.e. }}
\newcommand{\resp}{\emph{resp. }}
\newcommand{\vsp}{\vspace{0.5cm}}

\newcommand{\miniop}[3]{%
\renewcommand{\arraystretch}{0.6}
\begin{array}{c}
{\scriptstyle #1}\\
#2\\
{\scriptstyle #3}
\end{array}
\renewcommand{\arraystretch}{1}}

\newcommand{\Card}[1]{\vert #1 \vert}

\newcommand{\gap}{\Upsilon}
\newcommand{\1}{1\hspace{-1.3mm}1}

\newtheorem{theorem}{Theorem}[section]
\newtheorem{lemme}[theorem]{Lemma}
\newtheorem{defi}[theorem]{Definition}
\newtheorem{coro}[theorem]{Corollary}

\title[Chemical distance and first-passage percolation]{Asymptotic shape for the chemical distance and first-passage percolation in random environment}

\begin{document}

\authors{Olivier Garet,R{\'e}gine Marchand}
\address{Laboratoire de Math{\'e}matiques, Applications et Physique
Math{\'e}matique d'Orl{\'e}ans UMR 6628\\ Universit{\'e} d'Orl{\'e}ans\\ B.P.
6759\\
 45067 Orl{\'e}ans Cedex 2 France}
\email{Olivier.Garet@labomath.univ-orleans.fr}

\address{Institut Elie Cartan Nancy (math{\'e}matiques)\\
Universit{\'e} Henri Poincar{\'e} Nancy 1\\
Campus Scientifique, BP 239 \\
54506 Vandoeuvre-l{\`e}s-Nancy  Cedex France}
\email{Regine.Marchand@iecn.u-nancy.fr}


\begin{abstract}
The aim of this paper is to generalize the well-known asymptotic shape result for first-passage percolation on $\Zd$ to first-passage percolation on a random environment given by the infinite cluster of a supercritical Bernoulli percolation model.  We prove the convergence of the renormalized set of wet points to a deterministic shape that does not depend on the random environment.
As a special case of the previous result, we obtain an asymptotic shape theorem for the chemical distance in supercritical Bernoulli percolation.
We also prove a flat edge result. Some various examples are also given.
\end{abstract}

\maketitle

First-passage percolation was introduced in 1965 by Hammersley and Welsh~\cite{HW} as a model for the spread of a fluid in a porous medium. To each edge of the $\Zd$ lattice is attached a nonnegative random variable which corresponds to the travel time needed by the fluid to cross the edge.

When the passage times are independent identically distributed variables, Cox and Durrett showed that, under some moment conditions, the renormalized set of wet points at time $t$ almost surely converges to a deterministic asymptotic shape.
Deriennic~\cite{kesten}, and next Boivin~\cite{boivin}, progressively extended the result to the stationary ergodic case. 

In this paper, we want to study the analogous problem of spread of a fluid in a more complex medium.
On the one hand, an edge can either be open or closed according to the local properties of the medium -- e.g. according to the absence or the presence of non-porous particles.
In other words, the $\Zd$ lattice is replaced by a random environment given by the infinite cluster of a super-critical Bernoulli percolation model.
On the other hand, as in the classical model, a random passage time is attached to each open edge. This random time corresponds to the local porosity of the medium -- e.g. the density of the porous phase.
Thus, our model can be seen as a combination between classical Bernoulli percolation and stationary first-passage percolation.

Our goal is to prove in this context the convergence of the renormalized set of wet points to a deterministic shape that does not depend on the random environment.
As a special case of the previous result, we obtain an asymptotic shape theorem for the chemical distance in supercritical Bernoulli percolation.

In the first section,  we give some notations adapted to this problem, precise the assumptions on the passage time and introduce the background of ergodic theory needed for the main proof.
In section~2, we give some estimations on the chemical distance and the travel times.
In section~3, we prove the existence of an asymptotic speed in a given direction and study its properties ( e.g. continuity, homogeneity and sub-additivity\dots)
Section~4 is devoted to the question of the positivity of these speeds, and in section~5 we prove the asymptotic shape result. In section~6, we study the existence of a flat edge in the asymptotic shape. 
In the last section, we develop some examples, such as exponential times, correlated chi-square passage times, and a model for road networks.

\section{Notations, definitions and preliminary lemmas}

\subsection*{The percolation structure}

Let us first construct a Bernoulli percolation structure on
$\Z^d$. Consider the 
graph whose vertices are the points of 
$\Z^d$, and put a non-oriented edge between each pair $(x,y)$ of points in $\Z^d$ such
that the Euclidean distance between $x$ and $y$ is equal to $1$. Two such
points are called \textit{neighbours}. This set
of edges is denoted by $\Ed$. 

Set
$\Omega_E=\{0,1\}^{\Ed}$.
In the whole paper, $p$ is supposed to satisfy  to $$p\in (p_c,1],$$ where $p_c=p_c(d)$ is the critical threshold for bond
percolation on $\Z^d$.
We denote by $\P_p$ the product probability
$(p\delta_1+(1-p)\delta_0)^{\otimes \Ed}$. 

A point $\omega$ in $\Omega_E$
is a 
\textit{random environment} for the first passage percolation. An edge $e \in \Ed$ is said to be
\textit{open} in the environment $\omega$ if $\omega(e)=1$, and \textit{closed}
otherwise. The states of the different edges are thus
independent under $\P_p$ .

A \textit{path} is a sequence $\gamma=(x_1,
e_1,x_2,e_2,\ldots,x_n,e_n,x_{n+1})$ such that $x_i$ and $x_{i+1}$ are
neighbours and $e_i$ is the edge between $x_i$ and $x_{i+1}$. 
We will also sometimes describe $\gamma$ only by the vertices it
visits $\gamma=(x_1,x_2,\ldots,x_n,x_{n+1})$ or by its edges
$\gamma=(e_1,e_2,\ldots,e_n)$. The number $n$ of edges in $\gamma$ is called
the \textit{length} of $\gamma$ and is denoted by $|\gamma|$. Moreover, we
will only consider \textit{simple paths} for which the visited vertices are all
distinct. A path is said to be \textit{open} in the environment $\omega$
if all its edges are open in $\omega$.

The \textit{clusters} of a environment $\omega$ are the connected
components of the graph induced on $\Z^d$ by the open edges in
$\omega$. For $x$ in $\Z^d$, we denote by $C(x)$ the cluster containing
$x$. In other words, $C(x)$ is the set of points in $\Z^d$ that are linked
to $x$ by an open path. 
We note $x\communique y$ to signify that $x$ and $y$ belong to the same cluster.
For $p> p_c$, there exists almost surely one
and only one infinite cluster. We denote by $C_{\infty}$ the random set:
$C_{\infty}=\{k\in \Zd;\Card{C(k)}=+\infty\}$. Note that $C_\infty$ is
almost surely connected.

Since $p> p_c$, we have $\P_p(0\in C_\infty)>0$. In our first-passage percolation context, it is natural to work on the
infinite cluster. Thus we introduce the conditioned
probability measure
$\Pcond_p$ on $\Omega_E$ by 
$$\Pcond_p(A)=\frac{\P_p(A\cap \{0\in C_{\infty}\})}{\P_p(0\in
  C_{\infty})}.$$

\subsection*{The passage times}
We give now to each edge of the graph $\Z^d$ a \textit{passage time}, which
represents the time needed by a fluid to cross the edge. Consider thus
$\Omega_S=\R_{+}^{\Ed}$, and a probability measure $S_{\nu}$ on
$\Omega_S$. 
We use the abusive notation $S_{\nu}$ to help the reader to keep in mind
the case of a product measure $S_{\nu}=\nu^{\otimes\Ed}$, although we will
consider the weaker assumption that $S_{\nu}$ is a stationary ergodic probability measure -- see below for a precise definition.

\subsection*{The probability space}
Our probability space will then be $\Omega=\Omega_E\times \Omega_S$. 
A point in $\Omega$ will be
denoted $(\omega, \eta)$, with $\omega$ corresponding to the environment
and $\eta$ assigning to each edge a passage time. 
On $\Omega$, we consider
the probability $\P=\P_p\otimes S_{\nu}$
and its conditioned version:
$$\Pcond=\Pcond_p\otimes S_{\nu}.$$

\subsection*{The two distances}
We first introduce the \textit{chemical distance} $D(x,y)(\omega)$ between $x$ and
$y$ in $\Zd$, 
depending only on the Bernoulli percolation structure $\omega$:
$$D(x,y)(\omega)= \inf_{\gamma}|\gamma|,$$
where the infimum is taken on the set of paths whose extremities are $x$
and $y$ and that are open in the environment $\omega$. 
It is of course only defined when $x$ and $y$ are in the same percolation
cluster. Otherwise, we set by convention $D(x,y)=+\infty$.
The distance $D(x,y)$ is thus, when it is finite, the minimal number of
open edges needed to link $x$ and $y$ in the environment $\omega$. 
We denote by $\gamma(x,y)$ an open path
between $x$ and $y$ with this minimal number of edges. When necessary, we
can define uniquely $\gamma(x,y)$ by choosing on order on the set of edges
$\Ed$ and taking the minimal such path for the lexical order on the edges
of the path.

For $(\omega,\eta)\in\Omega$, and $(x,y)\in\Zd\times\Zd$, we define the
\textit{travel time}
$d(x,y)(\omega,\eta)$ to be 
$$\inf_{\gamma}\sum_{e\in \gamma}\eta_e,$$
where the infimum is taken on the set of paths whose extremities are $x$
and $y$ and that are open in the environment $\omega$. Of course $d(x,y)=+\infty$ if and only 
if $D(k,l)=+\infty$.
It follows that  the event $\{d(k,l)<+\infty\}$ only depends on the
Bernoulli percolation structure 
$\omega$.
For $t\ge 0$, we note
$$B_t=\{k\in\Zd; d(0,k)\le t\}.$$ It is the random set of points that can be
reached from the origin in a time smaller than $t$, and it is a
non-decreasing random set in
$t$.

\subsection*{Usual definitions and results in ergodic theory}
We  recall here some classical definitions, notations  and results in ergodic theory which
can be found in any course -- see for instance  the book
written by Brown \cite{Brown}. 

Consider a probability space $(\Omega,\mathcal{F},\Qjoli)$. A map $\theta:
\Omega \rightarrow \Omega $ is said to be a \textit{measure-preserving
  transformation} if for each $A\in\mathcal{F}$, $\Qjoli(\theta^{-1}(A))=\Qjoli(A)$.
The quadruple $(\Omega,\mathcal{F},\Qjoli,\theta)$ is then called a \textit{dynamical system}.

For any measurable set $A$, we can define the entrance time in $A$:
$$n_A(x)=\inf\{n\ge 1; \theta^n x\in A\}.$$
The well-known recurrence theorem of Poincar{\'e} ensures that
for $\Qjoli$ almost every $\omega\in A$,  $\{n\ge 1; \theta^n x\in A\}$ in infinite.

If every subset of $\Omega$ which is invariant by $\theta$ has probability
$0$ or $1$,
we say that $\theta$ is \textit{ergodic} for $\Qjoli$.
When $\theta$ is ergodic, the mean entrance time is given by Kac's formula:
$$\int_A n_A\ d\Qjoli=1.$$

For a set $A\in\mathcal{F}$ such that $\Qjoli(A)>0$, the entrance time in $A$
induces a transformation $\theta_A:\Omega\to\Omega$ defined by
$\theta_A(\omega)=\theta^{n_A(\omega)}(\omega)$. 
Of course, $\theta_A$ is in general only defined for almost every  $\omega \in A$ -- this is a consequence of the recurrence theorem --  but its definition may be
extended to $A$.
If $\Qjoli_A$ is defined by 
$\Qjoli_A(B)=\frac{\Qjoli(A\cap B)}{\Qjoli(A)}$, the following results hold:
\begin{itemize}
\item If $\theta$ is a  $\Qjoli$-preserving transformation, then $\theta_A$ is a  $\Qjoli_A$-preserving transformation.
\item If $\theta$ is ergodic for $\Qjoli$, then $\theta_A$ is ergodic for $\Qjoli_A$.
\end{itemize}
 Note that
Kac's formula can here be rewritten in the following form:
$$\int_{\Omega} n_A\ d\Qjoli_{A}=\frac1{\Qjoli(A)}.$$

\subsection*{The translation operators}
\begin{defi}
For any set $X$, and $u\in\Zd$, we define the translation operator $\theta_u$ 
on $X^{\Ed}$ by the relation
$$\forall e\in \Ed\quad (\theta_u\omega)_e=\omega_{u.e},$$
where $u.e$ denotes the natural action of $\Zd$ on $\Ed$: if $e=\{a,b\}$, then
$u.e=\{a+u,b+u\}$. 
\end{defi}

In the whole paper, the measure $S_{\nu}$ which describes the passage time is supposed to
be an ergodic measure preserving the $\Zd$ actions $(\theta_{u})_{u\in\Zd}$.
It means that for each $u\in\Zd$,
$(\Omega_S,\mathcal{B}(\Omega_S),S_{\nu},\theta_u)$ is an ergodic dynamical
system. This covers of course the case of a product measure $S_{\nu}=\nu^{\otimes\Ed}$.

Since $S_{\nu}$ is invariant under the action of $\Zd$, the law of the passage time of a given edge only depends from its
direction. 
Let us denote by $(e_1,\dots,e_d)$ the canonical basis of $\Zd$.
In the whole paper, we will suppose that

\vspace{0.2cm}
\begin{tabular}{p{1cm}p{20cm}}
$\displaystyle (H_{\text{int}})$ &  
$\displaystyle m=\miniop{}{\sup}{1\le i\le d}\int \eta_{(0,e_i)}\ dS_{\nu}(\eta)=\miniop{}{\sup}{e\in\Ed}\int \eta_{e}\ dS_{\nu}(\eta)<+\infty.$
\end{tabular} 
\vspace{0.2cm}

\noindent When needed, we will suppose that $S_{\nu}$ satisfies to 


\begin{tabular}{p{1cm}p{20cm}}
$\displaystyle (H_{\alpha})$ &  
$\displaystyle \exists A,B>0\ \text{such that }\quad\forall \Lambda\subseteq \Ed\quad S_{\nu}\left(\eta\in\Omega_S;\sum_{e\in\Lambda} \eta_i\ge B|\Lambda|\right)\le\frac{A}{| \Lambda|^{\alpha}}$
\end{tabular} 

\noindent for an appropriate value of $\alpha$. Assumption $(H_{\alpha})$ is in fact a combination between assumptions on the moments and assumptions on the dependence between the passage times of distinct edges.

For instance, if $S_{\nu}$ is the product measure $\nu^{\otimes\Ed}$,
assumption $(H_{\alpha})$  follows from the Marcinkiewicz-Zygmund inequality
as soon as the passage time of an edge has a moment of order $2\alpha$ -- see e.g. Theorem 3.7.8 in \cite{stout}.
Note that obviously, $(H_{\alpha})$ is always fulfilled when the passage times are bounded.

\begin{defi}
For each $u\in\Zd\backslash\{0\}$, let $T_u(\omega)=\inf\{n\ge 1; nu\in C_\infty(\omega)\}$ and define the
random translation operator on $\Omega=\Omega_E \times \Omega_S$
$$\Theta_u(\omega,\eta)=(\theta_u^{T_u(\omega)}(\omega),\theta_u^{T_u(\omega)}(\eta)).$$
\end{defi}
Note that $T_u$ only depends on the environment  $\omega$,
and not on the 
passage times $\eta$, whereas the operator $\Theta_u$ acts on the whole
configuration $(\omega,\eta)$.
The new configuration $\Theta_u(\omega,\eta)$ is the
initial configuration viewed from the first point of the form $nu$ to be in
the infinite cluster.

\begin{lemme} \label{lemme_invariance}
$\Theta_u$ is a $\Pcond$-preserving transformation.
Moreover, $\Theta_u$ is ergodic for $\Pcond$ and
$$\F T_u=\frac1{\P_p(0\in C_{\infty})}.$$
\end{lemme}

\begin{proof}
Let $u\in\Zd$.
The dynamical system $(\Omega,\mathcal{B}(\Omega),\P,\theta_u)$ is the direct product
of the mixing system $(\Omega_E,\mathcal{B}(\Omega_E),\P_p,\theta_u)$
by the ergodic system $(\Omega_S,\mathcal{B}(\Omega_S),S_{\nu},\theta_u)$.
Since the product of a weakly mixing system and an ergodic system is ergodic -- see for instance proposition 1.6 of Brown  \cite{Brown} -- , it follows
that  $\theta_u$ is ergodic for $\P$.
It is easy to see that $\Theta_u$ is the  transformation induced by $\theta_u$ on  the set $A=\{0\in C_{\infty}\}\times \Omega_S$.
Since $\P(A)=\P_p( 0\in C_{\infty})S_{\nu}(\Omega_S)=\P_p( 0\in C_{\infty})>0$,
the announced result follows from the classical results of ergodic theory.
\end{proof}

\begin{defi}
Let us define for $n\ge 1$
$$T_{n,u}(\omega)=\sum_{k=0}^{n-1} T_u(\Theta_u^k\omega).$$
\end{defi}
Remark that $T_{n,u}$ represents the position of the $n^{\text{th}}$ intersection of the
infinite cluster with the half line $\N^* u$. 

\begin{lemme}
\label{T_integrable}
$\displaystyle \frac{T_{n,u}}n\to\frac1{\P_p(0 \in
  C_\infty)}\quad\Pcond\as$ 
\end{lemme}
\begin{proof}
Since $\Theta_u$ is ergodic for $\Pcond$, this immediately follows from
the pointwise ergodic theorem.
\end{proof}

\section{Estimations on the chemical distance and travel times}

We first recall the fundamental result of Antal
and  Pisztora (Theorem 1.1 in \cite{AP96}):  If $p>p_c(d)$, there exist 
a constant $\rho=\rho(p,d)\in [1,+\infty)$ and two
strictly positive constants $A$ and $B$ such that:
\begin{equation} \label{equ_AP96}
\forall y \in \Zd, \; \P_p(0 \leftrightarrow y, \; D(0,y) > \rho \|y\|_1 )\leq A
\exp \left( -B\|y\|_1 \right).
\end{equation}
Note that if $y=(y_1,\ldots,y_d) \in \R^d$, $\|y\|_1$ is defined as follows:
$\|y\|_1=\sum_{i=1}^d |y_i|$.

\begin{lemme} \label{analogueAP}
If $(H_{\alpha})$ holds with two positive constants $A$ and $B$, and if $\rho$ is the constant given by (\ref{equ_AP96}), then
there exists a positive constant $A'$ such
that 
$$
\forall x \in \Z^d, \; \P(0\communique x;d(0,x)\ge B\rho\Vert x\Vert_1)\le \frac{A'}{\Vert x\Vert_1^{\alpha}}.$$
\end{lemme}

\begin{proof} Let $x \in \Z^d$.
By definition, $d(0,x) \leq \sum_{e \in
    \gamma(0,x)} \eta(e)$.
Then
$$\P(0\communique x;d(0,x)\ge B\rho\Vert x\Vert_1 \leq \P \left( 0\communique x;
\sum_{e \in \gamma(0,x)} \eta(e)\ge B\rho\Vert x\Vert_1 \right).$$
We can now cut this probability in two pieces:
\begin{eqnarray*}
&& \P \left(  0\communique x; \sum_{e \in \gamma(0,x)} \eta(e)\ge B\rho\Vert
  x\Vert_1 \right) \\
& \le & \P \left( 0\communique x; D(0,x)\le \rho\Vert x\Vert_1; \sum_{e \in \gamma(0,x)}
  \eta(e)\ge B\rho\Vert x\Vert_1 \right) 
 + \P \left( 0\communique x;D(0,x)> \rho\Vert x \Vert_1\right).
\end{eqnarray*}
Note that we already know that the second term decreases exponentially fast
thanks to (\ref{equ_AP96}); for the first one, we have:
\begin{eqnarray*}
&&
\P \left(0\communique x;   D(0,x)\le \rho\Vert x\Vert_1; \sum_{e \in \gamma(0,x)} \eta(e)\ge
 B\rho\Vert x\Vert_1 \right) \\
& = & \sum_{\gamma,|\gamma |\le \rho\Vert
  x\Vert_1}\P\left(0\communique x; \gamma(0,x)=\gamma; \sum_{e \in \gamma} \eta(e)\ge
  B\rho\Vert x\Vert_1 \right)\\
 & = & \sum_{\gamma,|\gamma |\le \rho\Vert x\Vert_1}\P_p
 (0\communique x; \gamma(0,x)=\gamma)S_\nu\left( \sum_{e \in \gamma} \eta(e)\ge B\rho\Vert
   x\Vert_1 \right)\\
 & \le & \sum_{\gamma,|\gamma|\le \rho\Vert x\Vert_1}\P_p (0\communique x; 
   \gamma(0,x)=\gamma)S_\nu\left( \sum_{e \in \gamma} \eta(e)\ge B|\gamma|
   \right)\\
 & \le & \sum_{\gamma,|\gamma |\le \rho\Vert
   x\Vert_1}\P_p(0\communique x; \gamma(0,x)=\gamma)\frac{A}{\Vert\gamma\Vert^\alpha} \; \; \mbox{
   with Hypothesis ($H_\alpha$)}\\
& \le & \frac{A\P_p(0\communique x)}{\Vert x\Vert_1^{\alpha}}  \; \; \mbox{
  as } |\gamma(0,x)| \geq \Vert x\Vert_1.\\
\end{eqnarray*}
\end{proof}

\begin{lemme} \label{estimee_petit}
There exists two positive constants $A',B'$, depending only on the dimension $d$, such that for every $r>0$, for every $x \in
\Z^d$ such that $\Vert x\Vert_1\le r$, 
$$
\P_p(0\communique x;D(0,x)\ge (3r)^d)  
 \le  A'\exp \left( -B'r \right).
$$
\end{lemme}

\begin{proof}
Take $\rho,A,B$ as given by (\ref{equ_AP96}).
Assume that 
$\displaystyle r \geq \max\left\{2, \left( \frac{10\rho}{3^d}
\right)^{\frac{1}{d-1}} \right\}$. Take $x \in \Z^d$ such that $\|x\|_1 \leq r$.  

We denote by $B(r)$ the cube $[-5r/4,5r/4]^d \cap \Z^d$. As the
segment 
$[a,b]$ contains at most $b-a+1$ integer points, 
$|B(r)| \leq (\frac{5r}{2}+1)^d,$ and, as $r \geq 2$,
$|B(r)|<(3r)^d$.
Denote by $\partial B(r)$ the external boundary of $B(r)$, defined as the set of
points in $\Zd \cap B(r)^c$ that admit a neighbour in $B(r)$. Note that $x
\in B(r)$ and that if $y \in \partial B(r)$, then $\|x-y\|_1
\geq \frac{r}{4}$.

Assume that $D(0,x)>(3r)^d$, then $D(0,x)>|B(r)|$, and then  $\gamma(0,x)$ must visit at least
one point in
$\partial B(r)$: 
\begin{eqnarray*}
& & \P(0 \leftrightarrow x, \;  D(0,x) \geq (3r)^{d}) \\
& \leq & \P \left( \exists y \in \partial B(r) \mbox{ such that } 0
  \leftrightarrow y, \; y \leftrightarrow x, \; D(0,y) 
  \geq \frac{(3r)^{d}}{2} \mbox{ or } D(x,y) \geq \frac{(3r)^{d}}{2} \right) \\
& \leq & \sum_{y \in \partial B(r)} \P \left( 0 \leftrightarrow y, D(0,y)
  \geq \frac{(3r)^d}{2} \right) + \P \left( x \leftrightarrow y, D(x,y)
  \geq \frac{(3r)^d}{2} \right) .
\end{eqnarray*} 
But $\displaystyle \frac{(3r)^d}{2} \geq 5\rho r \geq
\rho \max\{\|x\|_1,\|x-y\|_1\}$, and thus
\begin{eqnarray*}
& & \P(0 \leftrightarrow x, \;  D(0,x) \geq (3r)^{d}) \\& \leq & \sum_{y \in \partial B(r)} \P \left( 0 \leftrightarrow y, D(0,y)
  \geq \rho \|x\|_1 \right) + \P \left( x \leftrightarrow y, D(x,y)
  \geq \rho  \|x-y\|_1\right). 
\end{eqnarray*} 
But $\min \{\|x\|_1,\|x-y\|_1\} \geq
\frac{r}{4}$ and using (\ref{equ_AP96}) we get:
\begin{eqnarray*}
\P(0 \leftrightarrow x, \;  D(0,x) \geq (3r)^{d}) 
& \leq & 2 \sum_{y \in \partial B(r)} A \exp \left(-\frac{Br}{4} \right) \\
& \leq &  2 (2d)(5r/2+1)^{d-1} A \exp \left(-\frac{Br}{4}\right).
\end{eqnarray*} 
This ends the proof of the lemma.
\end{proof}

\begin{lemme} \label{estimee_petit_deux}
If $(H_{\alpha})$ holds, then
there exist positive constants $a$ and  $A'$, depending only of the
dimension $d$ of the grid,  such that for every $r>0$, for every $x \in
\Z^d$ satisfying $\|x\|_1 \leq r$, 
$$\P(0 \leftrightarrow x, d(0,x) \geq ar^d) \leq \frac{A'}{r^{\nu}},$$
with $\nu={\alpha-d+1}$.
\end{lemme}

\begin{proof}
Let $b>0$ be a constant that will be fixed later. Denote by
$B_r$ the cube $[-(1+b)r,(1+b)r]^d \cap Z^d$, and by $\partial_{ext} B_r$ the set of points in $\Z^d \backslash B_r$ that admit a neighbour inside $B_r$. 
We have 
$$|B_r| \leq (2(1+b)r+1)^d .$$
Denote by $n_e(B_r)$ the number of edges thqt have their two extremities in $B_r$. Then
$$d(2(1+b)r-2)^d \leq n_e(B_r) \leq d(2(1+b)r)^d.$$
Now, let $A$ and $B$ be the two constants given by Hypothesis $(H_\alpha)$, and choose $a>0$ such that 
$$\forall r>0, \; a^d > Bd(2(1+b)r)^d.$$
Take $x \in \Z^d$ such that $\|x\|_1 \leq r$.  Note that 
$x \in B_n$. Now,
\begin{eqnarray*}
& & \P  (0 \leftrightarrow x, \;  d(0,x) \geq ar^d) \\
&=& \P \left(0 \leftrightarrow x, \;  d(0,x) \geq ar^d \mbox{ and } \sum_{e \in B_r} \eta_e < Bd(2(1+b)r)^d \right) \\
& & + \P \left(0 \leftrightarrow x, \;  d(0,x) \geq ar^d \mbox{ and } \sum_{e \in B_r} \eta_e \geq Bd(2(1+b)r)^d \right)\\
& \leq & \P \left( \exists y \in \partial_{ext} B_n \mbox{ such that } 0
  \leftrightarrow y, \; y \leftrightarrow x, \; d(0,y) 
  \geq \frac{ar^d}{2} \mbox{ or } d(x,y) \geq \frac{ar^d}{2} \right) \\
&& + \P \left(\sum_{e \in B_r} \eta_e \geq Bn_e(B_r) \right).
\end{eqnarray*}
By Hypothesis $(H_\alpha)$, the last term is smaller than $\displaystyle \frac{A}{n_e(B_r)^\alpha}$. Note moreover that for
$r$ large enough:
$$B\rho\|y\|_1\leq B\rho((1+b)r +1) \leq  ar^d/2 \; \; \; \mbox{ and } \; \;
\;\|y\|_1 \geq br$$
$$B\rho\|x-y\|_1\leq B\rho((1+b)r +1 )\leq
ar^d/2 \; \; \; \mbox{ and } \; \;
\;\|y-x\|_1 \geq br.$$
Thus
\begin{eqnarray*}
&& \P \left( \exists y \in \partial_{ext} B_r \mbox{ such that } 0
  \leftrightarrow y, \; y \leftrightarrow x, \; d(0,y) 
  \geq \frac{ar^d}{2} \mbox{ or } d(x,y) \geq \frac{ar^d}{2} \right) \\
& \leq & \sum_{y \in \partial_{ext} B_r}  \P \left( 0 \leftrightarrow y, d(0,y)
  \geq \frac{ar^d}{2} \right) + \P\left( x \leftrightarrow y, d(x,y)
  \geq \frac{ar^d}{2} \right) \\
& \leq &\sum_{y \in \partial_{ext} B_r}  \P \left( 0 \leftrightarrow y, d(0,y)
  \geq  B\rho \|y\|_1\right) + \P \left( x \leftrightarrow y, d(x,y)
  \geq  B\rho \|x-y\|_1 \right) \\
& \leq & \sum_{y \in \partial_{ext} B_r}  A(\frac1{\|y\|_1^{\alpha}}+
\frac1{\|x-y\|_1^{\alpha}}) \;\;\; \mbox{ with Lemma \ref{analogueAP}}\\
 & \leq & 2d(2(1+b)r+1)^{{d-1}} \frac{A}{(br)^{\alpha}}.
\end{eqnarray*} 

\end{proof}

\begin{lemme} 
\label{integrabilite-et-psi}
There exist a function $\psi:\Zd\to \R^{+}$ satisfying $\lim_{\Vert
  u\Vert_1\to +\infty} \psi(u)=0$, and such that
$$\forall u\in\Zd,\quad \displaystyle \E_{\Pcond_p}(D(0,T_u u))\le \frac{\Vert u\Vert_1}{\P_p(0\communique\infty)}(\rho+\psi(u)).$$ 
\end{lemme}

\begin{proof}
Remember the following identity:
$$\E_{\Pcond_p}\frac{D(0,T_u u)}{\Vert u\Vert}  =   \int_{0}^{+\infty} \Pcond_p(D(0,T_u u)\geq h\Vert u\Vert) \ dh,$$
and  cut $\F(D(0,T_u u))$ in three terms, according to the values of $T_u$.
Let us estimate first the contribution for $T_u$ small:
\begin{eqnarray*}
S_1 (h,u) & = & \Pcond_p \left( 
D(0,T_u u) \geq h \Vert u\Vert_1, T_u \leq
  \frac{(h\Vert u\Vert_1)^{1/d}}{3 \Vert u\Vert_1} \right) \\
& = & \frac{1}{\P_p (0 \leftrightarrow \infty)} \P_p \left( 
D(0,T_u u) \geq h \Vert u\Vert_1, T_u \leq
  \frac{(h\Vert u\Vert_1)^{1/d}}{3 \Vert u\Vert_1} \right) \\
& = & \frac{1}{\P_p (0 \leftrightarrow \infty)} \sum_{ k \leq \frac{(h\Vert u\Vert_1)^{1/d}}{3 \Vert u\Vert_1} } 
\P_p
    \left(0 \leftrightarrow ku,  D(0,k u)\geq h\Vert u\Vert_1, T_u=k\right) \\
& \leq & \frac{1}{\P_p (0 \leftrightarrow \infty)} \sum_{k \leq \frac{(h\Vert u\Vert_1)^{1/d}}{3 \Vert u\Vert_1} } 
\P_p \left( 0 \leftrightarrow ku,  D(0,k u)\geq h\Vert u\Vert_1 \right) .
\end{eqnarray*}
But as $\displaystyle k \leq \frac{(h\Vert u\Vert_1)^{1/d}}{3 \Vert u\Vert_1}$, we have $h \Vert u\Vert_1\geq (3k \Vert u\Vert_1)^d$
and thus, by the previous lemma, 
\begin{eqnarray*}
S_1 (h,u)& \leq & \frac{1}{\P_p (0 \leftrightarrow \infty)} \sum_{k \leq \frac{(h\Vert u\Vert_1)^{1/d}}{3 \Vert u\Vert_1}} 
\P_p\left(0 \leftrightarrow ku,   D(0,k u)\geq (3k \Vert u\Vert_1)^d\right)\\
& \leq & \frac{1}{\P_p (0 \leftrightarrow \infty)} \frac{(h\Vert u\Vert_1)^{1/d}}{3 \Vert u\Vert_1} A' \exp \left( -B'\frac{(h\Vert u\Vert_1)^{1/d}}{3} \right).
\end{eqnarray*}
Now, $\int_{\R^+}S_1(h,u)dh$ is finite, and there exists a constant $C_1$,
independent from $u$, such that
\begin{eqnarray*}
\int_0^\infty S_1(h,u)dh 
& \leq &  \frac{1}{\P_p (0 \leftrightarrow \infty)} \int_0^\infty  \frac{(h\Vert
  u\Vert_1)^{1/d}}{3 \Vert u\Vert_1} A' \exp \left( -B'\frac{(h\Vert
    u\Vert_1)^{1/d}}{3} \right) dh \\
& \leq  & \frac{1}{\P_p (0 \leftrightarrow \infty)} \frac{A'}{3\Vert u\Vert_1^2}
\int_0^\infty x^{1/d} \exp\left( -B'\frac{x^{1/d}}{3} \right) dx \\
& \leq &  \frac{C_1}{\Vert u\Vert_1^2}.
\end{eqnarray*}
Let us now estimate the contribution for $T_u$ medium:
\begin{eqnarray*}
S_2 (h,u) & = & \Pcond_p \left( 
D(0,T_u u) \geq h \Vert u\Vert_1, \frac{(h\Vert u\Vert_1)^{1/d}}{3 \Vert
  u\Vert_1}< T_u \leq \frac{h }{\rho}\right) \\
& = & \frac{1}{\P_p (0 \leftrightarrow \infty)} \P_p \left( 
D(0,T_u u) \geq h \Vert u\Vert_1, \frac{(h\Vert u\Vert_1)^{1/d}}{3 \Vert
  u\Vert_1}< T_u \leq \frac{h }{\rho}\right) \\
 & = & \frac{1}{\P_p (0 \leftrightarrow \infty)} \sum_{ \frac{(h\Vert u\Vert_1)^{1/d}}{3 \Vert
  u\Vert_1}< k  \leq \frac{h}{\rho} } \P_p \left( 
D(0,T_u u) \geq h \Vert u\Vert_1, T_u=k \right) \\
& \leq & \frac{1}{\P_p (0 \leftrightarrow \infty)} \sum_{  \frac{(h\Vert u\Vert_1)^{1/d}}{3 \Vert
  u\Vert_1}< k  \leq \frac{h }{\rho} } \P_p \left( 
0 \leftrightarrow ku,  D(0,ku) \geq h \Vert u\Vert_1  \right) \\
& \leq & \frac{1}{\P_p (0 \leftrightarrow \infty)} \sum_{  \frac{(h\Vert u\Vert_1)^{1/d}}{3 \Vert
  u\Vert_1}< k  \leq \frac{h }{\rho} } \P_p \left( 
0 \leftrightarrow ku,  D(0,ku) \geq \rho k \Vert u\Vert_1 \right) \\
&& \; \; \; \; \; \mbox{ and with (\ref{equ_AP96}):}\\
& \leq & \frac{1}{\P_p (0 \leftrightarrow \infty)} \sum_{  \frac{(h\Vert u\Vert_1)^{1/d}}{3 \Vert
  u\Vert_1}< k  \leq \frac{h }{\rho} } A \exp \left( -B k \Vert u\Vert_1
\right) \\
& \leq & \frac{1}{\P_p (0 \leftrightarrow \infty)} \frac{h }{\rho} A \exp \left( -B \frac{(h\Vert u\Vert_1)^{1/d}}{3} \right). 
\end{eqnarray*}
Thus, there exists a constant $C_2$,
independent from $u$, such that
\begin{eqnarray*}
\int_0^\infty S_2(h,u)dh 
& \leq &  \frac{1}{\P_p (0 \leftrightarrow \infty)} \int_0^\infty  \frac{h }{\rho} A \exp \left( -B \frac{(h\Vert
    u\Vert_1)^{1/d}}{3} \right) dh \\
& \leq & \frac{1}{\P_p (0 \leftrightarrow \infty)} \frac{A}{\rho \|u\|_2} \int_0^\infty x^{1/d} \exp \left( -B
  \frac{x^{1/d}}{3} \right) dx \\
& \leq & \frac{C_2}{\Vert u\Vert_1^2}.
\end{eqnarray*}
Finally, for $T_u$ large:
\begin{eqnarray*}
S_3 (h,u) & = & \Pcond_p \left( 
D(0,T_u u) \geq h \Vert u\Vert_1, T_u > \frac{h }{\rho}\right) \\
& \leq &  \Pcond_p \left( T_u > \frac{h }{\rho}\right)
\end{eqnarray*}
and
\begin{eqnarray*}
\int_0^\infty S_3(h,u)dh & \leq &  \int_0^\infty\Pcond_p \left( T_u > \frac{h }{\rho}\right)dh \\
& \leq & \rho \E_{\Pcond_p} T_u \;\;\; \mbox{ with Lemma \ref{lemme_invariance}}\\
& \leq & \frac{\rho}{\P_p(0 \leftrightarrow \infty)}.
\end{eqnarray*}
Putting these three pieces together, we get the announced result.

\end{proof}

\begin{coro} \label{corollaire_integrable}
For every $u\in\Zd\backslash\{0\}$,
$$ \E_{\Pcond}(d(0,T_u u))\leq m \E_{\Pcond_p}(D(0,T_u u))<\infty.$$
\end{coro}
\begin{proof}
\begin{eqnarray*}
\E_{\Pcond}(d(0,T_u u)) 
& = & \sum_{k=0}^\infty \E_{\Pcond}(d(0,k) \1_{\{T_u=k\}}) \\
& \leq  & \sum_{k=0}^\infty \E_{\Pcond} \left(\1_{\{T_u=k\}} \sum_{e \in
    \gamma(0,ku)} \eta(e) \right) \\
& \leq  & \sum_{k=0}^\infty \E_{\Pcond} \left( \left. \1_{\{T_u=k\}} \E_{\Pcond} \left(
    \sum_{e \in \gamma(0,ku)} \eta(e) \right| \omega \right) \right) \\
& \leq  & \sum_{k=0}^\infty  \E_{\Pcond_p}\left( \1_{\{T_u=k\}} m
  D(0,ku) \right) \\
& \leq  & m \E_{\Pcond_p}(D(0,T_u u)).
\end{eqnarray*}
The result of the previous lemma gives then the integrability of $d(0,T_u u)$.
\end{proof}

\section{Existence of directional speed constants and related properties}

The next Lemma is the analogous of the first step in the classical study of first-passage percolation, that is the proof of the existence of an asymptotic speed in a given direction.

Unsurprisingly, the fundamental tool is here again Birkhof's subadditive ergodic theorem. The main difference is that it is applied to an induced dynamical system, and not directly. It is the reason why checking the integrability (lemma~\ref{corollaire_integrable}) was more intricate than in the classical case.

\begin{lemme}
For $n\ge 0$ and $u\in\Zd\backslash\{0\}$, we define the travel time $f_{n,u}$ between the origin and the
$n$-th intersection of the infinite cluster with $\Z_{+}^{*}u$:
$$f_{n,u}(\omega,\eta)=d(0,(T_{n,u}(\omega)u)(\omega,\eta).$$
Then there exists a constant $f_u\ge 0$ such that
$$\frac{f_{n,u}}n\to f_u\quad\Pcond\as$$
The convergence also holds in $L^1(\Pcond)$.
Moreover, $f_u\le \E_{\Pcond} d(0,T_u u)$.
\end{lemme}

\begin{proof} 
For the convenience of the reader, we now choose $u\in\Zd\backslash\{0\}$,
and we write $f_n$ (\resp $f$, $T_n$, $T$)
instead of $f_{n,u}$ (\resp $f_u$, $T_{n,u}$, $T_u$ ).

Note that the travel time is a sub-additive function, we thus have:
\begin{eqnarray*}
f_{n+k} & = & d(0,T_{n+k}u) \\
& \le & d(0,T_n u)+d(T_n u,T_{n+k} u) \\
& = & d(0,T_n u)+d(0,T_k u)\circ T_n \\
& = & f_n+f_k\circ T_n.
\end{eqnarray*}
Moreover, the $f_n$ are non-negative  and, thanks to lemma~\ref{lemme_invariance}, $\Pcond$ is invariant under $\Theta_u$.
Then, as  $f_1$ is integrable thanks to lemma~\ref{corollaire_integrable}, 
Birkhof's sub-additive ergodic theorem says there exists a $\Theta_u$-invariant function $f$ such that

$$\frac{f_n}n\to f \; \; \Pcond \as \mbox{ and in } L^1(\Pcond).$$
Since $\Theta_u$ is ergodic for $\Pcond$, every  $\Theta_u$-invariant function
is $\Pcond$ almost surely constant, so $f$ is constant.
Moreover, the deterministic sub-additive ergodic theorem ensures that 
$$\E_{\Pcond} f=\inf_{n\ge 0}\frac{\E_{\Pcond} f_n}n\le \E_{\Pcond} f_1=\E_{\Pcond} d(0,T_u).$$
\end{proof}

We can now define a function $\mu$ which is exactly the analogous of the function that appears in the classical case.

\begin{theorem} \label{th_dir_rat} 
Suppose that $(H_{\alpha})$ holds for some $\alpha>1$.
For each $u \in
\Zd \backslash \{0\}$, we set $\mu(u)=\P_p(0\communique\infty)f_u$, where $f_u$ is
  given by the previous lemma, and $\mu(0)=0$. 
Then 
$$\lim_{n \rightarrow \infty} \frac{d(0,T_{n,u}u)}{T_{n,u}}=\mu(u)
\quad\Pcond\as$$
Moreover, this function $\mu:\Zd \rightarrow [0, \infty)$ has the following properties:
\begin{enumerate}
\item
$\mu$ is symmetric: for any vector  $u \in \Zd$: $\mu(-u)=\mu(u)$.

\noindent If, moreover, $S_{\nu}$ is a product measure, then
$\mu$ is invariant under permutations or reflections of the coordinates in
$\Zd$.
\item
$\mu$ is homogeneous, in the following sense: for any integer $k \in \Z $
and any vector  $u \in \Zd$:
$$\mu(ku)=|k| \mu(u).$$
\item
$\mu$ is bounded, in the following sense: for any vector $u \in \Zd$:
$$\mu(u)\le \rho m\Vert u\Vert_1 \; \; \mbox{ and } \mu(u)\le \mu_{*} \Vert u\Vert_1,$$
with $\mu_{*}=\sup_{1\le i\le d} \mu(e_i)$.
\item
$\mu$ is sub-additive, in the following sense: 
for any 
vectors $u$ and $v$ in $\Zd$, 
$$\mu(u+v) \leq \mu(u)+\mu(v).$$
\item
$\mu$ is continuous, in the following sense:
for any 
vectors $u$ and $v$ in $\Zd$, 
$$|\mu(u) -\mu(v)| \leq \mu_{*} \|u-v\|_1.$$
\end{enumerate}
\end{theorem}

We delay the proof of this theorem to make some remarks and comments.

\vspace{0.3cm}
\noindent\textbf{Remarks.} 
\begin{itemize}
\item
We exactly recover the properties that the results of classical first-passage percolation let expect.
\item
Let us now discuss the meaning of the inequality $\mu(u)\le\rho m\Vert u\Vert_1$.
In classical first-passage percolation, we would have the inequality $\mu(u)\le m\Vert u\Vert_1$. Roughly speaking, we can say that the additional factor $\rho=\rho(p,d)$ corresponds to the bound given by corollary 1.3 of Antal and Pisztora \ref{equ_AP96} for the asymptotic ratio between the chemical distance and the $l_1$ distance on $\Zd$:
\begin{equation}
\label{corap}
\miniop{}{\limsup}{\Vert y\Vert_1\to +\infty}\frac{D(0,y)}{\Vert y\Vert_1}\1_{\{0\communique y\}}\le \rho(p,d)\ \P_p \as
\end{equation}
\item 
In the special case $S_{\nu}=(\delta_{1})^{\otimes\Ed}$, the travel times coincide with the chemical distances.
As a direct corollary of the previous theorem, we obtain the existence of asymptotic directed speeds for the chemical distance:
For each $y\in\Zd\backslash\{0\}$, there exists a constant $\mu(y)$ such that
\begin{equation*}
\lim_{\substack{n\to +\infty\\ 0\communique ny}}\frac{D(0,ny)}{n}=\mu(y)\ \P_p \as
\end{equation*}
Since $\mu(y)\le \rho(p,d)\Vert y\Vert_1$, this precises the equation (\ref{corap}) of Antal and Pisztora.

\end{itemize}

\begin{proof}

\vsp
\noindent
\underline{Step 1.}
The convergence result is clear since $T_{n,u}/n\to 1/{\P_p(0 \in C_\infty)}$ and $f_{n,u}/n\to
f_u\ge 0$. Consequently, we obtain the identity $\mu(u)=\P_p(0\communique\infty)f_u$. 

\vsp
\noindent
\underline{Step 2.} 
It is easy to see that
$\frac{d(0,nu)}n$ weakly converges to $$\P_p(0\communique\infty)\delta_{f_u}+(1-\P_p(0\communique\infty))\delta_{+\infty}.$$
By translation invariance $d(0,nu)=d(nu,0)$ has the same law that $d(0,-nu)$.
It follows that $f_{u}=f_{-u}$, and then that $\mu(u)=\mu(-u)$.

When $S_{\nu}$ is a product measure, Property (1) is obvious since the whole model is
invariant under permutations or reflections of the coordinates in
$\Zd$.

\vsp
\noindent
\underline{Step 3.} 
Let us prove now that for any integer $k \in \Z \backslash \{0\}$ and any
non null vector $u \in \Zd$, $\mu(ku)=|k| \mu(u)$. Note that 
$$\frac{f_{n,ku}}{n}=|k| \frac{f_{n,ku}}{n|k|};$$
but $({f_{n,ku}}/{n|k|})_{n\ge 1}$ is a subsequence of $({f_{n,\varepsilon
u}}/{n})_{n\ge 1}$ with $\varepsilon=k/|k|$. Thus
$$\lim_{n\to +\infty} \frac{f_{n,ku}}{n}=|k| \lim_{n\to
  +\infty}\frac{f_{n,\varepsilon u}}{n}=|k|\mu(\varepsilon u)=|k|\mu(u),$$
by the previous step.

\vsp
\noindent
\underline{Step 4.}
By the previous step $\mu(u)=\frac1{n}\mu(nu)$.
Then
\begin{eqnarray*}
\mu(u)& = & \Pcond(0\communique\infty)\frac{f_{nu}}{n}\\
 & \le & \frac{\E_{\Pcond}(d(0,T_{nu}nu)}n  \Pcond(0\communique\infty)\\
 & \le & \frac{m\E_{\Pcond}(D(0,T_{nu}nu)}n  \Pcond(0\communique\infty)\\
 & \le & \frac{m(\rho+\psi(nu)\Vert nu\Vert_1}n\text{ by lemma }(\ref{integrabilite-et-psi})\\
& \le & m(\rho+\psi(nu)\Vert u\Vert_1
\end{eqnarray*}
Since $\lim_{n\to +\infty} \psi(nu)=0$, it follows that
$\mu(u)\le m\rho \Vert u\Vert_1$. This proves the first part of (3).

\vsp
\noindent
\underline{Step 5.} Let us now prove the sub-additive property of $\mu$. Let
$u$ and $v$ be two fixed non null vectors in $\Zd$. The main difficulty
here, compared to the classical proof in standard first-passage percolation
(see Kesten) is that two given points in $\Z^d$ are not linked by an open
path in every configuration $\omega$. 

\begin{lemme}
Let $x$ and $y$ be two linearly independent vectors in $\Z^d$. Then
$$\frac{1}{n} \sum_{k=1}^n \1_{\{kx \leftrightarrow \infty\}}
\1_{\{ky \leftrightarrow \infty\}}  \rightarrow \P_p(0\leftrightarrow
\infty)^2 \; \; \mbox{almost surely under } \P_p.$$
\end{lemme}

\begin{proof}
Let $\alpha=\frac13\min(\Vert x\Vert_1,\Vert y\Vert_1,\Vert x-y\Vert_1)$
and define $Z_k=\1_{\{\alpha k\le C(kx)\}}\1_{\{\alpha k\le C(ky)\}}$.
\begin{eqnarray*}
\P_p(Z_k\ne \1_{\{kx \leftrightarrow \infty\}}
\1_{\{ky \leftrightarrow \infty\}}) & \le &  \P_p(\alpha k\le |C(kx)|<+\infty)+\P_p(\alpha k\le |C(ky)|<+\infty)\\ & \le & 2\P_p(\alpha k\le |C(0)|<+\infty).
\end{eqnarray*}
Since $p>p_c$, $C(0)\1_{\{|C(0)|<+\infty\}}$ is integrable, and we get
$$\sum_{k\ge 1}\P_p(Z_k\ne \1_{\{kx \leftrightarrow \infty\}}
\1_{\{ky \leftrightarrow \infty\}})<+\infty.$$ Hence, Borel-Cantelli
lemma ensures that, $\P_p$ almost surely, $Z_k$ coincides with  $\1_{\{kx \leftrightarrow \infty\}}
\1_{\{ky \leftrightarrow \infty\}}$ for every $k$ large enough.
It follows that
$\frac{1}{n} \sum_{k=1}^n \1_{\{kx \leftrightarrow \infty\}}
\1_{\{ky \leftrightarrow \infty\}}$
has the same behavior that
$\frac{1}{n} \sum_{k=1}^n Z_k$.
But $(Z_k)_{k\ge 1}$ is a sequence of independent uniformly bounded random variables. Then
$$\frac{1}{n} \sum_{k=1}^n \big( Z_k-\E Z_k\big)  \rightarrow 0 \; \; \P_p \as$$
Since $\alpha=\frac13\min(\Vert x\Vert_1,\Vert y\Vert_1,\Vert x-y\Vert_1)$,
 $$\E_{\P_p} Z_k=\E_{\P_p} \1_{\{\alpha k\le
   |C(kx)|\}}\E_{\P_p}\1_{\{\alpha k\le |C(ky)|\}}=\P_p(\alpha k\le
 C(0))^2\to \P_p(0\communique\infty)^2,$$ and the results follows from Cesaro's theorem.
\end{proof}

\begin{lemme} 
Suppose that there exists $\alpha>1$ such that $(H_{\alpha})$ is fulfilled
with two positive constants $A$ and $B$ and let $\rho$ be the constant
given by (\ref{equ_AP96}). 
Let $x$ and $y$ be two linearly independent vectors in
  $\Z^d$, and for $k \in \Z_+^{*}$, denote by
  $A_k$ the following event:
$$A_k=\{ d(0,kx) \leq {B\rho} k \|x\|_1, \; d(0,k(x+y)) \leq {B\rho} k \|x+y\|_1, \;
d(kx,k(x+y)) \leq {B\rho} k \|y\|_1 \}.$$
Then, almost
  surely under $\P$, if
\begin{eqnarray*}
u_n & = & \frac{1}{n} \sum_{k=1}^n \frac{d(0,kx)}{k} \1_{A_k}
\1_{\{0  \leftrightarrow \infty\}}
\1_{\{kx \leftrightarrow \infty\}}
\1_{\{k(x+y) \leftrightarrow \infty\}}, 
\end{eqnarray*}
then $(u_n)_n$ converges to $ \mu(x) \1_{\{0  \leftrightarrow \infty\}} \P_p
( 0  \leftrightarrow \infty)^2$ and moreover, $(\E_{\P_p}(u_n))_n$ converges to $\mu(x) \P_p
( 0  \leftrightarrow \infty)^3$.
\end{lemme}

\begin{proof}
Note first that
\begin{eqnarray}
&& \P_p \left( \1_{A_k}
\1_{\{0  \leftrightarrow \infty\}}
\1_{\{kx \leftrightarrow \infty\}}
\1_{\{k(x+y) \leftrightarrow \infty\}} 
\neq \1_{\{0  \leftrightarrow \infty\}}
\1_{\{kx \leftrightarrow \infty\}}
\1_{\{k(x+y) \leftrightarrow \infty\}} \right) \nonumber \\
& = & \E_{\P_p} \left( \1_{\{0  \leftrightarrow \infty\}}
\1_{\{kx \leftrightarrow \infty\}}
\1_{\{k(x+y) \leftrightarrow \infty\}} (1-\1_{A_k}) \right)
\nonumber \\
& \leq & \E_{\P_p} \left( \1_{\{0  \leftrightarrow \infty\}}
\1_{\{kx \leftrightarrow \infty\}}
\1_{\{k(x+y) \leftrightarrow \infty\}} \right.\nonumber \\
& & \times \left. \left( \1_{ \{d(0,kx)> {B\rho} k \|x\|_1\}}
+ \1_{ \{d(0,k(x+y))> {B\rho} k \|x+y\|_1\}}
+ \1_{ \{d(kx,k(x+y))> {B\rho} k \|y\|_1\}}\right)
\right) \nonumber \\
& \leq & 
\P_p \left( d(0,kx)> {B\rho} k \|x\|_1 \mbox{ and } 0  \leftrightarrow kx
\right) 
\nonumber \\
&  & +\P_p \left(d(0,k(x+y))> {B\rho} k \|x+y\|_1 \mbox{ and } 0
  \leftrightarrow k(x+y) \right)  \nonumber \\
& & +\P_p \left( d(kx,k(x+y))> {B\rho} k \|y\|_1 \mbox{ and } kx
  \leftrightarrow k(x+y) \right)  \nonumber \\
& \leq & \frac{A}{k^{\alpha}} \left( \frac{1}{\|x\|_1^{\alpha}} + \frac{1}{\|x+y \|_1^{\alpha}}
  +  \frac{1}{\|y\|_1^{\alpha}} \right) \; \; \;  \mbox{
   with the previous lemma.} \nonumber 
\end{eqnarray}
Then, using Borel-Cantelli Lemma, we know
that almost surely there exists $K=K(\omega)>0$ such that for any $k \geq K$,
$$\1_{A_k}\1_{\{0  \leftrightarrow \infty\}}
\1_{\{kx \leftrightarrow \infty\}}
\1_{\{k(x+y) \leftrightarrow \infty\}} = \1_{\{0  \leftrightarrow \infty\}}
\1_{\{kx \leftrightarrow \infty\}}
\1_{\{k(x+y) \leftrightarrow \infty\}}.$$ 

Choose $\varepsilon>0$. Then, as a consequence of the convergence result in Theorem~\ref{th_dir_rat}, and
enlarging $K=K(\omega)$ if necessary, one has, almost surely:
$$\forall k \geq K(\omega), \; \left| \frac{d(0,kx)}{k}-\mu(x) \right| \1_{\{0  \leftrightarrow \infty\}}
\1_{\{kx \leftrightarrow \infty\}} \leq \varepsilon.$$
Now, for every $n$:
\begin{eqnarray}
&& \left| u_n -\mu(x) \1_{\{0  \leftrightarrow \infty\}} \P_p
( 0  \leftrightarrow \infty)^2 \right| \nonumber \\
& \leq & \frac{1}{n}  \sum_{k=1}^{K(\omega)-1} \frac{d(0,kx)}{k} 
\1_{A_k}
\1_{\{0  \leftrightarrow \infty\}}
\1_{\{kx \leftrightarrow \infty\}} 
\1_{\{k(x+y) \leftrightarrow \infty\}} \nonumber \\
& + & \frac{1}{n}  \sum_{k=K(\omega)}^n \left| \frac{d(0,kx)}{k}-\mu(x)
\right| 
\1_{A_k}
\1_{\{0  \leftrightarrow \infty\}}
\1_{\{kx \leftrightarrow \infty\}} 
\1_{\{k(x+y) \leftrightarrow \infty\}} \nonumber \\
& + & \mu(x) \left| \frac{1}{n}  \sum_{k=K(\omega)}^n \1_{A_k}
\1_{\{0  \leftrightarrow \infty\}}
\1_{\{kx \leftrightarrow \infty\}} 
\1_{\{k(x+y) \leftrightarrow \infty\}} 
-\1_{\{0  \leftrightarrow \infty\}} \P_p
( 0  \leftrightarrow \infty)^2 \right| .\nonumber 
\end{eqnarray}
But we can choose $N(\omega)$ large enough to ensure that if $n \geq N(\omega)$, then the first term in the sum is smaller than
$\varepsilon$. The second one is smaller than $\varepsilon$ thanks to the
choice of $K(\omega)$. For the last term, remember that the choice of
$K(\omega)$ ensures that for every $k \geq K(\omega)$,
we can forget the indicator function $\1_{A_k}$. But then the previous lemma ensures the almost sure
convergence to $0$ of the last term.
This proves the desired almost-sure convergence.

To prove the convergence in mean, note that $u_n \leq \rho$ and use the
dominated convergence theorem.
\end{proof}

Let us now prove the sub-additive property of $\mu$. Note first that if $x$
and $y$ are linearly dependent, the inequality is in fact an equality
thanks the the homogeneity property. Now, let $x$ and $y$ be two linearly
independent vectors in $\Z^d$. For every $k \geq 1$, let $A_k$ be the same
event as the one 
defined in the previous lemma. Then, for every $k \geq 1$,
\begin{eqnarray}
&& \frac{d(0,k(x+y))}{k} \1_{A_k}
\1_{\{0  \leftrightarrow \infty\}}
\1_{\{kx \leftrightarrow \infty\}}
\1_{\{k(x+y) \leftrightarrow \infty\}}  \nonumber \\
& \leq & \frac{d(0,kx)}{k} \1_{A_k}
\1_{\{0  \leftrightarrow \infty\}}
\1_{\{kx \leftrightarrow \infty\}}
\1_{\{k(x+y) \leftrightarrow \infty\}} \nonumber \\
& + & \frac{d(kx,k(x+y))}{k} \1_{A_k}
\1_{\{0  \leftrightarrow \infty\}}
\1_{\{kx \leftrightarrow \infty\}}
\1_{\{k(x+y) \leftrightarrow \infty\}} \nonumber 
\end{eqnarray}
Summing for $1 \leq k \leq n $ and taking the means gives:
\begin{eqnarray}
&& \E_{\P_p} \left( \frac{1}{n} \sum_{k=1}^n 
\frac{d(0,k(x+y))}{k} \1_{A_k}
\1_{\{0  \leftrightarrow \infty\}}
\1_{\{kx \leftrightarrow \infty\}}
\1_{\{k(x+y) \leftrightarrow \infty\}}  \right) \label{a} \\
& \leq & \E_{\P_p} \left( \frac{1}{n} \sum_{k=1}^n 
\frac{d(0,kx)}{k} \1_{A_k}
\1_{\{0  \leftrightarrow \infty\}}
\1_{\{kx \leftrightarrow \infty\}}
\1_{\{k(x+y) \leftrightarrow \infty\}} \right) \label{b} \\
& + & \E_{\P_p} \left( \frac{1}{n} \sum_{k=1}^n 
\frac{d(kx,k(x+y))}{k} \1_{A_k}
\1_{\{0  \leftrightarrow \infty\}}
\1_{\{kx \leftrightarrow \infty\}}
\1_{\{k(x+y) \leftrightarrow \infty\}} \right) \label{c} 
\end{eqnarray}
Thanks to the previous lemma, (\ref{a}) and (\ref{b}) converge respectively
to $\mu(x+y)\P_p
( 0  \leftrightarrow \infty)^3$ and $\mu(x)\P_p
( 0  \leftrightarrow \infty)^3$. Moreover, 
\begin{eqnarray}
\mbox{(\ref{c})} & = & \frac{1}{n} \sum_{k=1}^n \E_{\P_p} \left(
  \frac{d(kx,k(x+y))}{k} \1_{A_k} 
\1_{\{0  \leftrightarrow \infty\}}
\1_{\{kx \leftrightarrow \infty\}}
\1_{\{k(x+y) \leftrightarrow \infty\}} \right) \nonumber \\
& = & \frac{1}{n} \sum_{k=1}^n \E_{\P_p} \left( 
\frac{d(0,ky)}{k} \1_{B_k}
\1_{\{0  \leftrightarrow \infty\}}
\1_{\{-ky \leftrightarrow \infty\}}
\1_{\{-k(x+y) \leftrightarrow \infty\}} \right), 
 \nonumber 
\end{eqnarray}
where $B_k$ is the event
$$B_k=\{ d(0,-ky) \leq \rho B k \|y\|, \; d(0,-k(x+y)) \leq \rho B k \|x+y\|, \;
d(-ky,-k(x+y)) \leq \rho B  k \|x\| \}.$$ This was obtained by translating
along the vector $-kx$. Once again, the previous lemma enables us to
conclude that (\ref{c}) converges to $\mu(y)\P_p
( 0  \leftrightarrow \infty)^3$. As $\P_p
( 0  \leftrightarrow \infty)>0$, this leads to the desired inequality.

\vsp
\noindent
\underline{Step 6.} Property (5) is a direct consequence of the
sub-additivity of $\mu$.

\vsp
\noindent
\underline{Step 7.} Using properties (1), (2) and (4), we get:
$$\mu(u_1,\ldots,u_d) \leq \sum_{i=1}^d |u_i| \mu(e_i)$$
where $(e_i)_{1 \leq i \leq d}$ is the canonical basis of $\Zd$. Thus, 
$$\mu(u) \leq \mu_{*} \|u\|_1.$$
The second part of property (3) follows immediately.\end{proof}

This finally allows us to prove an analogue result for any rational
direction. Indeed, let $q \in \Q^d
\backslash \{0\}$ be fixed. Choose any couple $(N,u) \in \N^* \times \Z^d$
such that $u=Nq \in
\Zd$ and set
$\mu(q)=\mu(u)/N$. The constant $\mu(q)$ is well defined
thanks to the homogeneity property of $\mu$. The function $\mu$, extended in this manner to
$\Q^d \backslash \{0\}$, obviously keeps its properties: 
for any rational $\alpha$ and any
vectors $u$ and $v$ in $\Q^d$, 
$$\mu(\alpha u)=|\alpha| \mu(u), \; \mbox{ and } \mu(u+v) \leq
\mu(u)+\mu(v) \mbox{ and } |\mu(u) -\mu(v)| \leq \rho m\|u-v\|_1.$$
Thus we can extend $\mu$ to $\R^d$ by continuity, keeping these three
properties.

\section{Sufficient conditions for the positivity of $\mu$}

As in the classical cas, we study now the question of the positivity of $\mu$, which is a crucial point in the quest of an asymptotic shape result.
We begin with the case of a product measure.

\begin{theorem}
\label{norme-indep}
Suppose here that $S_{\nu}=\nu^{\Ed}$.
Then
\begin{itemize}
\item If  $p \nu(0) <p_c$, then $\mu$ is a norm on $\R^d$.
\item If $p \nu(0) >  p_c$, then $\mu=0$ on $\R^d$.
\end{itemize}
\end{theorem}

\begin{proof}
The only remaining point to consider is to determine whether $\mu(u)=0$
implies $u=0$ or not. 
Let us define, for each $e\in\Ed$,  $\eta'_e=\eta_e\omega_e+(1-\omega_e)$.
The law of $(\eta'_e)_{e\in\Ed}$ under $\P$ is
$(p\nu+(1-p)\delta_1)^{\otimes\Ed}$. We can define the new distance
relatively to these new passage times:
$$d'(x,y)(\omega,\eta)=\inf_{\gamma}\sum_{e\in \gamma}\eta'_e,
$$ 
where the infinimum is taken over all the paths from $x$ to $y$ in $\Zd$.
This is in fact the classical
 first-passage percolation distance associated to the passage times $\nu'=
 p\nu+(1-p)\delta_1$. Denote by $\mu'$ its associated renormalized
 limit. Note that by construction $d'(x,y)\le d(x,y)$, and thus that $\mu'
 \leq \mu$. It
 is a now classical result (see Kesten \cite{kesten} for instance) that
 $\nu'(0)<p_c$ implies that the associated limit $\mu'$ is a norm on
 $\R^d$, and, since $\nu'(0)=p\nu(0)$, if $p\nu(0)<p_c$ then $\mu$ is a norm on
 $\R^d$.

On the other hand, suppose that $p\nu(0)>p_c$. Color in red the edges $e$
such that $\omega(e)=1$ and $\nu(e)=0$. As $p\nu(0)>p_c$, there is almost
surely a unique infinite cluster of red edges, included in the infinite
cluster for the  Bernoulli$(p)$ percolation structure. Let $u\in\Zd$. On the event $R=$ ``the
origin is in the red infinite cluster'' (included in the event ``the origin
is in the Bernoulli$(p)$ infinite cluster''), 
we can find an increasing sequence of (random) integers $(k_n)$ such that $k_n.u$ is in the red
infinite cluster for every $n$. But then clearly, $d(0, k_n.u)=0$. Thus,
evaluating $\mu(u)$ on this subsequence gives $\mu({u})=0$ on the event $R$. 
Since $\mu(u)$ is constant under $\Pcond$, this completes the proof.
\end{proof}

\begin{coro}
The asymptotic speed $\mu$ associated to the chemical distance in Bernoulli percolation with parameter $p>p_c(d)$ is a norm on $\Rd$. 
\end{coro}
When $S_{\nu}$ is not a product measure, there is no general method to determinate whether $\nu$ is a norm or not, even in the classical case $p=1$.

A natural idea is to use stochastic comparison: if $S_{\nu'}\succ S_{\nu}$,
then $\mu(p,\nu')\ge \mu(p,\nu)$.
Consequently, if $a$ is such that $S_{\nu'}(\eta_e\in [a,+\infty))=1$ holds for each $e\in\Ed$, then $S_{\nu'}$ is a norm for each $p>p_c$, because $d(.,.)\ge aD(.,.)$.

We will now prove that $\mu$ is a norm if $S_{\nu}$ satisfies  an appropriate large deviation inequality.

\begin{theorem}
\label{une-norme}
Let us suppose that there exists $A>0$ and a map $f:(0,+\infty)\to\R^{+}$ 
\begin{equation*}
\forall \epsilon>0\quad\forall \Lambda\subseteq \Ed\quad S_{\nu}\left(\eta\in\Omega_S;\sum_{e\in\Lambda} \eta_i\le \epsilon|\Lambda|\right)\le A\exp(-f(\epsilon)| \Lambda|). 
\end{equation*}
Let $K_0=\miniop{}{\liminf}{\epsilon\to 0^{+}} \exp(-f(\epsilon)).$ 

Then, the application $\mu$ associated to $S_{\nu}$ is a norm as soon as  
$K_0 p \lambda(d)<1$,
where $\lambda(d)$ is the $d$-dimensional connective constant
\end{theorem} 

\begin{proof}
Let  $\epsilon>0$ be small enough to ensure that $(\lambda(d)+\epsilon)\exp(-f(\epsilon))<1$
and consider $u\in\Zd$ .

\begin{eqnarray*}
\P(d(0,u)\le\epsilon\Vert u\Vert_1) & =  & \sum_{\gamma } \P(\gamma(0,u)=\gamma;
\sum_{e\in\Lambda} \eta_i\le  \epsilon\Vert u\Vert_1)\\
 & \le  & \sum_{\gamma } \P(\gamma(0,u)=\gamma;
\sum_{e\in\Lambda} \eta_i\le  \epsilon |\gamma |)\\
& \le & \sum_{n=\Vert u\Vert_1}^{+\infty} H (\lambda(d)+\epsilon)^{n}\exp(-nf(\epsilon))\\
& \le & \frac{H}{1-((\lambda(d)+\epsilon)\exp(-f(\epsilon))}\big((\lambda(d)+\epsilon)\exp(-f(\epsilon)) \big)^{\Vert u\Vert_1}.
\end{eqnarray*}
Then, it follows from Borel-Cantelli lemma that
$d(0,nu)\ge n\Vert u\Vert_1$ if $n$ is large enough: we can conclude
that $\mu(u)\ge\epsilon\Vert u\Vert_1$.
As $\mu$ is homogeneous and continuous, this inequality is extended first to
$\Q^d$, and next to $\Rd$.
Since $\Vert{.}\Vert_1$ is a norm, so does $\mu$.

\end{proof}

\noindent\textbf{Remarks.}
We will give in subsection~\ref{chi2-et-cie} an example of application of this theorem.

Let us now discuss the assumption of theorem~\ref{une-norme} by considering the case $S_{\nu}=(\Ber (1-q))^{\otimes \Ed}$.
By Chernof's theorem, one can take $$f(\epsilon)=(1-\epsilon)\log\frac{1-\epsilon}q+ \epsilon\log\frac{\epsilon}{1-q},$$ thus $K_0=q$.
Then, theorem~\ref{une-norme} says that $\mu$ is a norm as soon as $qp<\frac1{\lambda(d)}$ whereas theorem~\ref{norme-indep} teaches us that the optimal condition
is $qp<p_c$.\\ Since $\frac1{\lambda(d)}\sim \frac1{2d}\sim p_c(d)$ for large $d$, 
assumption of theorem~\ref{une-norme} does not seem so bad.

\section{The asymptotic shape theorem}

From now on, we suppose that $\mu$ is a norm. We want now study the convergence of the renormalized set of wet points at time $t$ to the unit ball for $\mu$.
Because of the presence of holes in the infinite cluster, the Hausdorff distance seems naturaly adapted to study the convergence.

\begin{defi}
For $x\in\Rd$ and $r\ge 0$, we define
$$B_{\mu}(x,r)=\{y\in\Rd;\mu(x-y)\le r\}.$$
\end{defi}

\begin{defi}
The Hausdorff distance between two non empty compact subsets  $A$ and $B$ of $\Rd$
is defined by
$$\mathcal{D}(A,B)=\inf\{r\ge 0; A\subset B+B_{\mu}(0,r)\text{ and } B\subset A+B_{\mu}(0,r)\}.$$
\end{defi}

Note that we use the Hausdorff distance associated to $\mu$, but the equivalence of norms on $\Rd$ ensures that the induced topology does not depend on this choice.

\begin{theorem}[asymptotic shape]
\label{asymptotic-shape}
Suppose that $(H_{\alpha})$ holds for some $\alpha>d^2+2d-1$ and that
$\mu$ is a norm.
 Then, 
$$\lim_{t\to +\infty}\mathcal{D}\left(\frac{B_t}t,B_{\mu}(0,1)\right)=0\quad\Pcond\as$$
\end{theorem}

\noindent\textbf{Remark.} It is important to have in mind that $(H_{\alpha})$ is always satisfied if the passage times are bounded by an absolute constant. We have already seen that $S_{\nu}$ is always a norm when the passage times are bounded from below by a a positive constant. Putting these two results together, we get
that $S_{\nu}$ always satisfies  the assumptions of theorem~\ref{asymptotic-shape} if it is an invariant ergodic measure for which we can found $0<a<b<+\infty$ with $S_{\nu}([a,b]^{\Ed})=1$. We thus obtain an asymptotic shape theorem for the chemical distance in supercritical Bernoulli percolation. Remember that on the event ``the origin is in the infinite cluster'', we denote by $B_n$ the set of points in $\Z^d$ whose chemical distance from $0$ is less or equal to $n$.
\begin{coro}
For every $p>p_c(d)$, there exists a deterministic convex compact set $A$ with non-empty interior and invariant under permutations or reflections of the coordinates in $\R^d$ such that:
$$\frac{B_n}{n} \to A \text{ for the Hausdorff topology } \Pcond_p \as$$
\end{coro}

We will first give some preliminary lemmas in order to prove theorem~\ref{asymptotic-shape}. 

\begin{lemme} \label{lemme_existence_de_points}
Let $z \in \Z^d$ be a fixed non null vector. Then for every $\varepsilon>0$,
$$ \Pcond (\exists N>0, \; \forall x\in [N,+\infty), \; \exists k \in [(1-\varepsilon)x,\ldots,(1+\varepsilon)x]\cap\N \; \mbox{ such that } kz \leftrightarrow \infty)=1.$$
\end{lemme}

\begin{proof}
Almost surely under $\Pcond$, the sequence $(T_{m,z})_{m\ge 1}$ is an unbounded increasing sequence; thus there exists $m\in\N$,
with  $T_{m,z}\le x<T_{m+1,z}$.
Now, 
$$\frac{x-T_{m,z}}{x}\le \frac{T_{m+1,z}-T_{m,z}}{T_{m,z}}= \frac{T_{m+1,z}}{T_{m,z}}-1=\frac{m+1}{m}\frac{ \frac{T_{m+1,z}}{m+1}}{ \frac{T_{m,z}}{m}}-1,$$
which vanishes at infinity since  $\frac{T_{m,z}}{m}$ tends to $\frac1{\P(0\communique\infty)}$.
It follows that $\frac{x-T_{m,z}}{x}\le\epsilon$ as soon as $x$ is large enough, thus we can take $k=T_{m,z}$.
\end{proof}

\begin{lemme} \label{lemme_distance_pas_trop_grande}
Suppose that $(H_{\alpha})$ holds with $\alpha>d^2+2d-1$ and
let $\varepsilon>0$ be small enough. Then 
$$\P(\exists M>0,\; \forall m \geq M, \; \forall x \in \Z^d, \forall y \in \Z^d,\;  (\|x\|_1=m,\|y\|_1 \leq \varepsilon m )\Rightarrow d(x,y) \leq B\rho \varepsilon m)=1.$$
\end{lemme}

\begin{proof}
In view to use Borel-Cantelli lemma, denote by $A_m$ the set:
$$A_m=\{ \exists x \in \Z^d, \; \exists y \in \Z^d  \mbox{ such that } \|x\|_1=m,\|y\|_1 \leq \varepsilon m \mbox{ and } d(x,y) \geq B\rho \varepsilon m\}.$$
Using lemma~\ref{estimee_petit_deux} and lemma~\ref{analogueAP}, we get (note that in the following lines, $K_i$ and $L_j$  denote constants, depending only on the dimension $d$, whose values are not precised to avoid intricate expressions).
We set $r=\left(\frac{B\rho\epsilon m}a\right)^{1/d}$.
\begin{eqnarray*}
\P(A_m) & \leq & \sum_{x, \|x\|_1=m} \sum_{l=0}^{\varepsilon m} \sum_{y, \|x-y\|_1=l} \P(x \leftrightarrow y, d(x,y) \geq B\rho \varepsilon m) \\
& \leq & \sum_{x, \|x\|_1=m} \sum_{l=0}^{r} \sum_{y, \|y\|_1=l} \P(y \leftrightarrow 0, d(y,0) \geq B\rho\epsilon m) \\
&& + \sum_{x, \|x\|_1=m} \sum_{l={r}}^{\varepsilon m} \sum_{y, \|y\|_1=l} \P(y \leftrightarrow 0, d(0,y) \geq B\rho \|y\|_1) \\
& \leq & \sum_{x, \|x\|_1=m} \sum_{l=0}^{r}  \sum_{y, \|y\|_1=l}  \frac{K_0}{r^{\nu}} \\
&& + \sum_{x, \|x\|_1=m} \sum_{l={r}}^{\varepsilon m} \sum_{y, \|y\|_1=l} \frac{L_0}{l^{\alpha}} \\
\end{eqnarray*}

Let us compute the first term
\begin{eqnarray*}
 \sum_{x, \|x\|_1=m} \sum_{l=0}^{r}  \sum_{y, \|y\|_1=l}  \frac{K_0}{r^{\nu}} & \le &
K_1 \sum_{x, \|x\|_1=m} \sum_{l=0}^{r}  l^{d-1}  \frac{1}{r^{\nu}} 
\\ & \le & K_2 \sum_{x, \|x\|_1=m}  \frac{r^d}{r^{\nu}} \\
 & \le & K_3 \frac{1}{m^{\frac{\nu}d-d}},
 \end{eqnarray*}

and the second term:
\begin{eqnarray*}
\sum_{x, \|x\|_1=m} \sum_{l={r}}^{\varepsilon m} \sum_{y, \|y\|_1=l} \frac{L_0}{l^{\alpha}} & \le & \sum_{x, \|x\|_1=m} \sum_{l={r}}^{\varepsilon m} l^{d-1} \frac{L_1}{l^{\alpha}}  \\
 & \le & \sum_{x, \|x\|_1=m} \sum_{l={r}}^{+\infty} l^{d-1} \frac{L_1}{l^{\alpha}}  \\
 & \le & \sum_{x, \|x\|_1=m} \frac{L_2}{m^{(\alpha-d)/d}}\\
 & \le &  \frac{L_3}{m^{\frac{\alpha-d}d-d+1}}.
\end{eqnarray*}

Remember that $\nu=\alpha-d+1$. The choice we made on $\alpha$ ensures that 
both terms are the general terms of  convergent series, which allows us to conclude thanks to Borel-Cantelli lemma.
\end{proof}

\begin{lemme}
\label{convergence_uniforme}
Let $\varepsilon>0$. Then 
$$\Pcond(\exists M>0, \; \forall y \in \Z^d, \; (\|y\|_1 \geq M\text{ and }y\communique 0 \Rightarrow |d(0,y)-\mu(y)| \leq \varepsilon \|y\|_1))=1.$$
\end{lemme}

\begin{proof}
Let us suppose this is false and choose $\varepsilon>0$ such that the assertion fails. Then, with positive probability under $\Pcond$, there exists a random sequence $(y_n)_{n \ge 0}$ of points in $\Z^d$ such that 
$$
\left\lbrace
\begin{array}{l} 
\|y_n\|_1 \rightarrow \infty \\
y_n \leftrightarrow 0 \\
|d(0,y_n)-\mu(y_n)| \geq \varepsilon \|y_n\|_1
\end{array}
\right.
$$
By considering a subsequence if necessary, suppose that
$$\frac{y_n}{\|y_n\|_1} \rightarrow z.$$
Let us approximate $z$ by a renormalized integer vector: consider  $\varepsilon_1>0 $ (to be chosen small enough later), and choose $z' \in\Zd$ be such that $\| \frac{z'}{\|z'\|_1} -z\|_1\leq \varepsilon_1$. 

Let us find, for each $y_n$, an integer point on the line $\R z'$ close enough from $y_n$:  let $h_n$ be the integer part of $\frac{\|y_n\|_1}{\|z'\|_1}$. We have
\begin{eqnarray*} 
\|y_n-h_n.z'\|_1 & \leq & \left\|y_n- \frac{\|y_n\|_1}{\|z'\|_1} z'\right\|_1+ \left| \frac{\|y_n\|_1}{\|z'\|_1}- h_n\right| \|z'\|_1 \\
& \leq &\|y_n\|_1  \left\|\frac{y_n}{\|y_n\|_1}-\frac{z'}{\|z'\|_1}\right\|_1 +\|z'\|_1
\end{eqnarray*}
Choose $N>0$ large enough to ensure that $n \geq N \Rightarrow \|\frac{y_n}{\|y_n\|_1}-z\|_1\leq \varepsilon_1$, \\
With the choice we wade for $z'$, we have 
$n \geq N \Rightarrow \left\|\frac{y_n}{\|y_n\|_1}-\frac{z'}{\|z'\|_1}\right\|_1 \leq 2\varepsilon_1$, whence
$$ \|y_n-h_n.z'\|_1 \leq 2\varepsilon_1\|y_n\|_1 +\|z'\|_1.$$

Lemma~\ref{lemme_existence_de_points} ensures now that, enlarging $N$ if necessary, there exists for every $n \geq N$, a integer $k_n \in [(1-\varepsilon_1)h_n,\ldots, (1+\varepsilon_1)h_n]$ such that $kz'$ is in the infinite cluster. But:
\begin{eqnarray*} 
\|y_n-k_n.z'\|_1 & \leq & \|y_n-h_n.z'\|_1 + |h_n-k_n| \|z'\|_1 \\
& \leq &  2\varepsilon_1\|y_n\|_1 +\|z'\|_1 +\varepsilon_1  h_n \|z'\|_1 \\
& \leq &  2\varepsilon_1\|y_n\|_1 +\|z'\|_1 +\varepsilon_1\|y_n\|_1 \\
& \leq &  3\varepsilon_1\|y_n\|_1 +\|z'\|_1 \\
& \leq &  4\varepsilon_1\|y_n\|_1,
\end{eqnarray*}
if $N$ is large enough. Enlarging once again $N$ if necessary, we can use lemma~\ref{lemme_distance_pas_trop_grande} to have: for every $n \geq N$: 
$$d(y_n,k_n.z') \leq 4 B\rho \varepsilon_1\|y_n\|_1.$$ 


Finally, for every $n$ large enough, we have:
\begin{eqnarray*} 
|d(0,y_n)-\mu(y_n)| & \leq & |d(0,y_n)-d(0,k_n.z')|+|d(0,k_n.z')-\mu(k_n.z')|+ 
|\mu(k_n.z')-\mu(y_n)| \\
& \leq & 4B \rho \varepsilon_1\|y_n\|_1 + k_n \left|\frac{d(0,k_n.z')}{k_n}-\mu(z')\right| +  \mu_1\|k_n.z'-y_n\|_1 \\
& \leq & 4B \rho \varepsilon_1\|y_n\|_1 +(1+\varepsilon_1)\frac{\|y_n\|_1}{\|z'\|_1}\left|\frac{d(0,k_n.z')}{k_n}-\mu(z')\right|+ 4\varepsilon_1\mu_1\|y_n\|_1.\end{eqnarray*}
But the convergence in the direction given by $z'$ ensures that for every $n$ large enough, 
$$\left|\frac{d(0,k_n.z')}{k_n}-\mu(z')\right| \leq \varepsilon_1.$$
So by choosing $\varepsilon_1$ small enough, we can ensure that for every $n$ large enough, 
$$|d(0,y_n)-\mu(y_n)| \leq \varepsilon\|y_n\|_1,$$
which yields to a contradiction.
\end{proof}

\begin{proof}[Proof of theorem~\ref{asymptotic-shape}]$\;$
Let $\epsilon>0$.

\noindent\underline{Step 1.}
If
$$ \exists T>0\quad \forall t \geq T\quad  \frac{B_t}{t} \subset B_{\mu}\left(0,(1+\varepsilon)\right)$$
fails, then there exists an unbounded increasing  sequence of times $(t_n)_{n\ge 1}$ and a sequence of vertices
$(y_n)_{n\ge 1}$
with $y_n\in {B(t_n)}$ and $\frac{y_n}{t_n}\notin B_{\mu}(0,1+\epsilon)$,
which gives $d(0,y_n)\le t_n$ and $\mu(y_n)>(1+\epsilon)t_n$.
By lemma~\ref{convergence_uniforme}, this happens with a null probability.

\noindent\underline{Step 2.}
Let us show that with probability 1, the following property holds:
\begin{equation}
\label{pas_de_trou}
 \exists T>0\quad \forall t \geq T\quad  B_{\mu}\left(0,1-\frac{\epsilon}2\right)\subset \frac{B_t}{t}+B_{\mu}\left(0,\frac{\varepsilon}2\right).
\end{equation}
If (\ref{pas_de_trou}) fails, 
then, with positive probability under $\Pcond$, there exists  an unbounded increasing  sequence of random times $(t_n)_{n\ge 1}$ and a sequence $(v_n)_{n\ge 1}$ of points in $\R^d$ such that 
$$
\left\lbrace
\begin{array}{l} 
v_n \in B_{\mu}(0,1-\frac{\epsilon}2)\\
v_n\notin \frac{B_{t_n}}{t_n}+B_{\mu}(0,\frac{\varepsilon}2)\\
\end{array}
\right.
$$
Let $v$ be a limiting value for $(v_n)_{n\ge 1}$: clearly, there exists an
unbounded increasing  sequence of times $t'_n$ such that
$v\notin \frac{B_{t'_n}}{t'_n}+B_{\mu}(0,\frac{\varepsilon}3)$.
Let $w\in\Q^d$ such that $\mu(v-w)\le \frac{\epsilon}{12}$.
Then, for each $n\ge 1$, we have
$w\notin \frac{B_{t'_n}}{t'_n}+B_{\mu}(0,\frac{\varepsilon}4)$, whence
$w\in B_{\mu}(0,1-\frac{5}{12}\epsilon)$.
We can write $w$ as $w=z/Z$, with $z\in\Zd$ and $Z\in\N$.

Lemma~\ref{lemme_existence_de_points} ensures
that, provided that $n$ is large enough, there exists
$z_n\in B_{\mu}(t'_n w,\frac{\epsilon}{4}\mu(t'_n w))$ with $z_n\communique\infty$. 
The vertex $z_n$ can not be in  $\frac{B_{t'_n}}{t'_n}$ because $w\notin \frac{B_{t'_n}}{t'_n}+B_{\mu}(0,\frac{\varepsilon}4)$. It follows that
$d(0,z_n)>t_n$.
But $$\mu(z_n)\le\mu(t'_n w)+\frac{\epsilon}{4}\mu(t'_n w)\le t'_n\mu(w)+\frac{\epsilon}{4}t'_n\le t'_n\left(1-\frac{5\epsilon}{12}+\frac{\epsilon}{4}\right)=\left(1-\frac{\epsilon}6\right)t'_n.$$
By lemma~\ref{convergence_uniforme}, this can (almost surely) not happen.
Then, (\ref{pas_de_trou}) holds almost surely.
Since
$$ B_{\mu}\left(0,1-\frac{\epsilon}2\right)\subset \frac{B_t}{t}+B_{\mu}\left(0,\frac{\varepsilon}2\right) \Longrightarrow  B_{\mu}(0,1)\subset \frac{B_t}{t}+B_{\mu}(0,{\varepsilon}),$$
this completes the proof.
\end{proof}

\section{Properties of the asymptotic shape}
As in the classical case, one has immediate properties for the asymptotic shape (see \cite{kesten}, chapter 6):

\begin{lemme}[Stochastic comparison]
\label{compsto}
 Let $p_c<p\leq p'$ and $\nu,\nu'$ be two probability measures on $\R^+$ such that the assumptions of the shape theorem are fulfilled for $S_\nu=\nu^{\otimes \Ed}$ and $S_{\nu'}=(\nu')^{\otimes \Ed}$ and such that $\nu'$ is stochastically smaller than $\nu$, then if $A$ and $A'$ denote their respective associated asymptotic shapes, we have
$$A \subset A'.$$
\end{lemme}

\begin{proof}
Coupling argument.
\end{proof}

\begin{lemme} Let $p>p_c$, and $\nu$ be a probability measure on $\R^+$ satisfiying the assumptions of the shape theorem. Denote by $\nu_{\min}$ the essential infimum of $\nu$,  $\nu_{mean}$ the mean of $\nu$, by $\mu_{p,\nu}$ the norm associated to the first passage-percolation model with $S_\nu=\nu^{\otimes \Ed}$ in the random environment given by the infinite cluster of Bernoulli$(p)$ percolation, and by $\tilde{\mu}_p$ the norm associated to the chemical distance for Bernoulli$(p)$ percolation. Then:
$$ \forall y \in \Z^d \; \; \; \nu_{\min} \tilde{\mu}_p (y)\leq \mu_{p,\nu} (y)\leq  \nu_{mean} \tilde{\mu}_p(y).$$
\end{lemme}

\begin{proof}
The lower bound is obtained by a coupling argument, the upper one by estimating the travel times on the minimal paths for the chemical distance.
\end{proof}

\begin{theorem}[Flat edge result]
\label{flatedge}
Here, we work on $\Z^2$. Let $p>p_c=1/2$. Denote by $\overrightarrow{p_c}$ the critical threshold for oriented percolation on $\Z^2$. 
Suppose that $S_\nu=\nu^{\otimes \Edeux}$, where $\nu$ is a probability measure on $\R^+$ such that $\nu_{\min}$, the essential infimum of $\nu$, is positive, satisfies the assumption of the shape theorem. Denote by $A_p$ the associated asymptotic shape, and by $\Diamond (r)$ the ball with radius $r$ for $\|.\|_1$.
\begin{itemize}
\item
$\displaystyle A_p \subset \Diamond  \left( \frac{1}{\nu_{\min}}\right).$
\item
If $p \nu(\nu_{\min})<\overrightarrow{p_c}$, then $\displaystyle A_p \subset \text{ int } \left(\Diamond  \left( \frac{1}{\nu_{\min}}\right)\right).$
\item
If $p \nu(\nu_{\min})>\overrightarrow{p_c}$, then 
$$
\displaystyle A_p \cap \text{ fr } \left(\Diamond  \left( \frac{1}{\nu_{\min}}\right)\right) \cap (\R^+)^2 =\left[\frac{1}{\nu_{\min}}M_q,\frac{1}{\nu_{\min}}N_q\right],$$
where 
$$\left\{
\begin{array}{l}
q=p \nu(\nu_{\min}),\\
\alpha_q \mbox{ is the asymptotic speed in supercritical oriented percolation } \\
\; \; \; \; \mbox{ with pamameter } q, \\
M_q=\left( \frac12+\frac{\alpha_q}{\sqrt{2}}, \frac12-\frac{\alpha_q}{\sqrt{2}} \right), \\
N_q=\left( \frac12-\frac{\alpha_q}{\sqrt{2}}, \frac12+\frac{\alpha_q}{\sqrt{2}} \right), \\
\left[\frac{1}{\nu_{\min}}M_q,\frac{1}{\nu_{\min}}N_q\right] \mbox{ denotes the segment line with extremities } \frac{1}{\nu_{\min}}M_q,\frac{1}{\nu_{\min}}N_q.
\end{array}
\right.
$$
\end{itemize}
\end{theorem}

\begin{coro}[Flat edge for the asymptotic shape for the chemical distance]
Here, we work on $\Z^2$. Let $p>p_c=1/2$. Denote by $\overrightarrow{p_c}$ the critical threshold for oriented percolation on $\Z^2$. 
Denote by $A_p$ the associated asymptotic shape, and by $\Diamond$ the unit ball with radius $1$ for $\|.\|_1$.
\begin{itemize}
\item
$\displaystyle A_p \subset \Diamond.$
\item
If $p <\overrightarrow{p_c}$, then $\displaystyle A_p \subset \text{ int } \left(\Diamond\right).$
\item
If $p >\overrightarrow{p_c}$, then 
$$
\displaystyle A_p \cap \text{fr} \left(\Diamond \right) \cap (\R^+)^2 =\left[M_p,N_p\right],$$
with the same notations as in the previous theorem.
\end{itemize}
\end{coro}

\newpage

The following pictures represent the set $B_n$ at time $n=450$ for the chemical distance. Note the absence of flat edge for $p=0.55<\overrightarrow{p_c}$.

\begin{tabular}{ccc}
\includegraphics[scale=1]{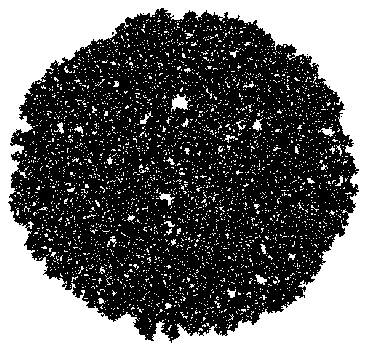} &  &
\includegraphics[scale=1]{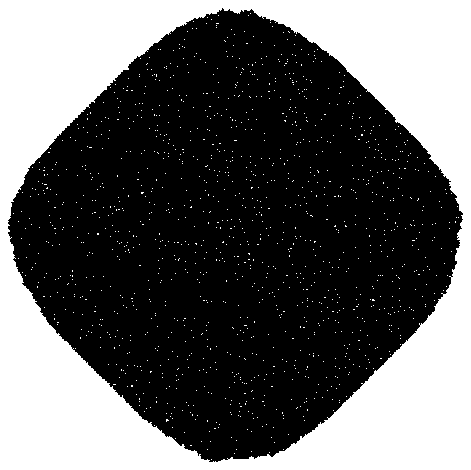}\\
$p=0.55$ & \hspace{-2cm} & $p=0.7$
\end{tabular}

\begin{proof}[Proof of theorem~\ref{flatedge}]
Note that the first point is a direct consequence of  lemma \ref{compsto}, by comparing our model with the classical case and passage times with law 
$\delta_{\nu_{\min}}$.

Note then, by lemma \ref{compsto}, that the asymptotic shape of the first-passage percolation for some $\nu^{\otimes \Edeux}$ on the infinite cluster of Bernoulli$(p)$ percolation is included in the asymptotic shape of classical first-passage percolation for the same $\nu^{\otimes \Edeux}$.  Thus the second point is a consequence of theorem 1.3 (ii) in \cite{marchand}, which is a consequence of the strict inequalities for the time constant in classical first-passage percolation established by van den Berg and Kesten in \cite{vdb-kes}.

For the third point,  let us examinate first the proof of the inclusion
$$A_p \cap \text{ fr } \left(D \left( \frac{1}{\nu_{\min}}\right)\right) \cap (\R^+)^2  \supset \left[\frac{1}{\nu_{\min}}M_q,\frac{1}{\nu_{\min}}N_q\right].$$
As noted in \cite{kesten} in the proof of theorem 6.13, the proof of the analogous result in the special case of Richardson's model treated in \cite{dur-lig} lemma 6-13 by Durrett and Liggett can immediately be adapted. 

It remains now to prove the exact length of the flat edge by proving:
$$A_p \cap \text{ fr } \left(D \left( \frac{1}{\nu_{\min}}\right)\right) \cap (\R^+)^2  \subset \left[\frac{1}{\nu_{\min}}M_q,\frac{1}{\nu_{\min}}N_q\right].$$
But once again, the stochastic comparison lemma and the exact determination result in the classical case -- theorem 1.3 in \cite{marchand} -- give the announced result.
\end{proof}

\section{Examples}

\subsection{Exponential passage times}
This subsection is devoted to the special case where the passage times
are independent variables with an exponential distribution $\mathcal{E}(\lambda)$. When $p=1$, we recover a Richardson model, which is one of the first studied examples for first-passage percolation -- see \cite{richardson}.

By theorem~\ref{norme-indep}, $\mu$ is a norm for $S_{\nu}=(\mathcal{E}(\lambda))^{\otimes\Ed}$ and each $p>p_c$.
Let us denote by $A(p,\lambda)$ the asymptotic form corresponding to
$S_{\nu}=(\mathcal{E}(\lambda))^{\otimes\Ed}$.
By the scaling properties of the exponential laws, it is easy to see
that $A(p,\lambda)=\lambda A(p,1)$ holds for each $\lambda>0$ and $p>p_c$. 
Let us prove that $A(p,1)\subset A(1,p)=p A(1,1)$.

The main tool is an appropriate coupling.
Let $(\eta_e)_{e\in\Ed}$, $(x_e)_{e\in\Ed}$, $(z_e)_{e\in\Ed}$
be independent random variables with for each $e\in\Ed$: $\eta_e\sim \mathcal{E}(p)$, $x_e\sim \mathcal{E}(1-p)$ and $z_e\sim \mathcal{E}(1)$.
Now define
$$\eta'_e=\1_{\{\eta_e\le x_e\}}\eta_e+\1_{\{\eta_e> x_e\}}z_e$$
and $$\omega_e=\1_{\{\eta_e\le x_e\}}.$$
On one hand,
\begin{eqnarray*}
P(\eta'_{e}>t,\omega_e=1) &  = & P(\eta_e>t,\eta_e\le x_e)\\
                       &  = & \int_{t}^{+\infty} pe^{-py}P(x_e>y)\ d\lambda(y)\\
 & = & \int_{t}^{+\infty} pe^{-py} e^{-(1-p)y}\ d\lambda(y)\\
 & = & pe^{-t}.
\end{eqnarray*}
Particularly, $P(\omega_e=1)=P(\eta'_e>t,\omega_e=1)=p$.
On the other hand
$$P(\eta'_e>t,\omega_e=0)= P(z_e>t,\eta_e\le x_e)=P(z_e>t)P(\eta_e\le x_e)=e^{-t}(1-t)$$
Thus, $\eta'_e$ and $\omega_e$ are independent random variables, with
$\eta'_e\sim\mathcal{E}(1)$ and $\omega_e\sim\Ber(p)$.

Consider now the set of wet points $B_t$ and $B'_t$ respectively associated to
$((1,\eta_{e})_{e\in\Ed})$ and $((\omega_e,\eta'_e))_{e\in\Ed}$.
Since $\eta'_e=\eta_e$ when the bond $e$ is open, it follows that
$$B'(t)\subset B(t)$$ always holds.\\ Dividing by $t$, and letting $t$ tend to
infinity, we get
$A(p,1)\subset A(1,p)=p A(1,1)$.

\subsection{Dependent $\chi^2$ passage times}
\label{chi2-et-cie}
We then give an example of unbounded dependent passage times which lead to
an asymptotic shape.
It is also an example in which our theorem~\ref{une-norme} makes us able to prove that $\mu$ is effectively a norm.

\begin{theorem}
Let $(X_n^i)_{n\in\Zd,i\in\{1,\dots,d\}}$ be a centered Gaussian process
such that 
$$\forall i,j\in\{1,\dots,d\}\quad\forall k,n\in\Zd\quad \E X_n^i X_{n+k}^j=\E X_0^i X_k^j$$
and
$$\sum_{k\in\Zd} \sum_{1\le i\le j\le d} |\E X_0^i X_k^j|<+\infty.$$
To avoid degenerate cases, we also suppose that
$$\forall i\in\{1,\dots,d\}\quad \E (X^i_0)^2>0.$$
Now consider the passage time $(\eta_{e})_{e\in\Ed}$ defined by
$\eta_{\{k,k+e_i\}}=(X_k^i)^2$.
If we denote by $S_{\nu}$ the law of $(\eta_{e})_{e\in\Ed}$, then 
$S_{\nu}$ satisfies to the assumption of theorem~\ref{asymptotic-shape}
for each $p>p_c$, so there is convergence to the asymptotic shape
for the associated growth model.
\end{theorem}

\begin{proof}
It is immediate that $S_{\nu}$ is shift-invariant. Since the covariance is summable, it tends to $0$
at infinity; it follows that $S_{\nu}$ is mixing, hence it is ergodic -- this result is attributed to Maruyama and Fomin by Sinaï \cite{sinai}.

Let us first recall some useful estimates whose proofs can be found in \cite{garet}.

\begin{lemme}
\label{borne_norme}
Let $X$ be a $\R^n$-valued centered Gaussian vector with covariance matrix $C$.
Let us denote by $\rho (C)$ the spectral radius \ie the largest norm of an eigenvalue of $C$.
Then, for each  $a^2>\rho(C)$, we have:
\begin{equation}
P(\Vert X\Vert_1^2\ge n a^2)\le 
e^{-nh(\frac{a^2}{\rho(C)})},
\end{equation}
where $h(x)=\frac12(x-\ln x-1)$.
The function $h$ is increasing and positive on $(1,+\infty)$,
with  $+\infty$ as limit at $+\infty$.
\end{lemme}

\begin{lemme}
\label{trouspectral}
Let $X$ be a $n$-dimensional  Gaussian  vector with positive definite covariance matrix $C$. Let us denote by $\gap (C)$ the spectral gap \ie the smallest eigenvalue of $C$. Then, for each
$a^2<\gap (C)$, we have
\begin{equation}
P(\Vert X\Vert_1^2\le n a^2)\le 
e^{-nh(\frac{a^2}{\gap (C)})},
\end{equation}
where $h(x)=\frac12(x-\ln x-1)$.
The function $h$ is positive and decreasing on $(0,1)$, with an infinite limit at 0.
\end{lemme}

We will now prove that $(\eta_e)_{e\in\Ed}$ satisfy to the assumption of 
 theorem~\ref{une-norme}.

Let us for convenience  identify $\Ed$ with $\Zd\times\{1,\dots, d\}$: the point $(x,i)\in\Zd\times\{1,\dots d\}$ will be identified with
$\{x,x+e_i\}\in\Ed$.
We define $$S=\sum_{k\in\Zd} \sum_{1\le i\le j\le d} |\E X_0^i X_k^j|<+\infty\text{ and }\sigma^2= \inf\{\E (X^i_0)^2>0;i\in\{1,\dots,d\}\}.$$
Let $\Lambda\subseteq \Zd\times\{1,\dots d\}$ and $\epsilon>0$.
We must bound
$$S_{\nu}\left(\eta\in\Omega_S;\sum_{e\in\Lambda} \eta_e\le \epsilon|\Lambda|\right)=\P\left(\sum_{(k,i)\in\Lambda} (X_k^i)^2\le \epsilon|\Lambda|\right).$$
Since $\displaystyle \sum_{k\in N\Zd\backslash\{0\}} \sum_{1\le i\le j\le d} |\E X_0^i
X_k^j| \le \sum_{k\ge N} \sum_{1\le i\le j\le d} |\E X_0^i X_k^j|$, we can find $N\in\N$ such that
$$\sum_{k\in N\Zd\backslash\{0\}} \sum_{1\le i\le j\le d} |\E X_0^i X_k^j| \le \frac{\sigma^2}2.$$
For each $(k,i)\in\{0,N-1\}^d\times \{1,\dots,d\}$, we can define
$A_{k,i}=\Lambda\cap \big((k+N\Zd)\times\{i\}\big)$.
By the pigeon-hole principle, there exists $(k_0,i_0)$ such that $\Card{A_{k_0,i_0}}\ge\frac{\Card{\Lambda}}{d N^d}$. 
Let $\tilde{X}$ be the $\Card{A_{k_0,i_0}}$-dimensional Gaussian vector composed by the
$(X_e)_{e\in A_{k_0,i_0}}$; it is obvious that
\begin{eqnarray*}
\P\left(\sum_{(k,i)\in\Lambda} (X_k^i)^2\le |\Lambda|\epsilon\right)
& \le & \P(\Vert \tilde{X}\Vert_2^2\le |\Lambda|\epsilon)\\
& \le & \P(\Vert \tilde{X}\Vert_2^2\le |A_{k_0,i_0}|{d N^d}\epsilon).\\
\end{eqnarray*}
By lemma~\ref{trouspectral},
$${d N^d}\epsilon<\gap (C) \Longrightarrow \P(\Vert \tilde{X}\Vert_2^2\le |A_{k_0,i_0}|{d N^d}\epsilon)\le\exp\left(-\Card{A_{k_0,i_0}}h\left(\frac{{d N^d}\epsilon}{\gap (C)}\right)\right),$$
where the covariance matrix for $\tilde{X}$ is $C=(\E X_k^iX_l^j)_{((k,i),(l,j))\in A_{k_0,i_0}\times A_{k_0,i_0}}$.
But
\begin{eqnarray*}
\gap (C) & \ge & \inf_{(k,i_0)\in A_{k_0,i_0}} \E (X_k^{i_0})^2-\sum_{(l,i_0)\in A_{k_0,i_0};(l,i_0)\ne (k,i_0)}\vert\E X_k^i X_l^{i_0} \vert \\
& \ge & \sigma^2-\sup_{(k,{i_0})\in A_{k_0,i_0}}\sum_{(l,{i_0})\in A_{k_0,i_0};(l,{i_0})\ne (k,i)} \vert\E  X_k^{i_0} X_l^{i_0}\vert\\
& \ge & \sigma^2-\sup_{(k,{i_0})\in A_{k_0,i_0}}\sum_{(l,{i_0})\in A_{k_0,i_0};(l,{i_0})\ne (k,i)} \vert\E  X_0^{i_0} X_{k-l}^{i_0}\vert\\
& \ge & \sigma^2-\sum_{k\in N\Zd\backslash\{0\}}\vert 
\E  X_0^{i_0} X_{k}^{i_0}\vert\\
& \ge & \sigma^2-\frac{\sigma^2}2=\frac{\sigma^2}2 .
\end{eqnarray*}
It follows that
\begin{equation}
\forall \epsilon\in \left(0,\frac{\sigma^2}{2Nd}\right)\quad P\left(\sum_{(k,i)\in\Lambda} (X_k^i)^2\le |\Lambda|\epsilon \right)\le\exp\left(
-\frac1{dN^d} h(\frac{2\epsilon}{\sigma^2})
\Card{\Lambda}
\right).
\end{equation}

We can  now claim that the assumptions of  theorem~\ref{une-norme} are fulfilled. Since $\lim_{x\to 0} h(x)=+\infty$, we have $K_0=0$.
Then, the application $\mu$ associated to $p$ and $S_{\nu}$ is a norm for each $p>p_c$.

With the help of lemma~\ref{borne_norme}, it is not difficult to see that
$$S_{\nu}\left(\eta\in\Omega_S;\sum_{e\in\Lambda} \eta_e\ge |\Lambda| a^2\right)\le
e^{-|\Lambda| h(\frac{a^2}{S})}$$
holds for each $a>S$.
It follows that $(H_{\alpha})$ holds for each $\alpha>1$, and hence that
the assumptions of theorem~\ref{asymptotic-shape} are fulfilled.
\end{proof}

\subsection{Cost ans speed for road networks.}

Imagine that there exist $n$ companies $(C_i)_{1 \le i \le n}$ that want to build roads on the same $\Z^d$ lattice: Company $C_i$, independentlty from the others, tries to build a road on edge $e$ with probability $p_i$ (set then $\omega_{e,i}=1$) and leaves the edge $e$ empty with probability $1-p_i$ (set then $\omega_{e,i}=0$), all constructions being independent. 
When two  or more companies want to build a road on the same edge, they share the possession of the road. 
A client of company $C_i$ can only use roads at least partially own by $C_i$. 
Suppose that each $p_i>p_c(d)$ to ensure that almost surely, company $C_i$ has a real (infinite) road network. Portions of roads own by comp $C_i$ that are not linked to the infinite cluster cannot be used by its clients. Set, if $e=\{a,b\}\in\Ed$ is an edge:
$$\eta_{e}=\sum_{i=1}^n \omega_{e,i}\1_{|C_i(a)|+|C_i(b)|=+\infty}.$$
Then $\eta_e$ counts the number of companies that offer the edge $e$ as a possible road for their clients.
Now, associate to the edge $e$ a random value:
$$\eta_{e,i}=f(i,\eta_e),$$
where $f$ is a deterministic function.

For instance, we can imagine that $ f(i,\eta_e)$ represents the costs for company $C_i$ to build a road on the edge $e$: this cost depend on the specificity of company $C_i$ (that is why $f$ depends on $i$), and also on the total number $\eta_e$ of companies that want to build a raod on th edge $e$ (their can share the costs of construction for instance). In this case, the function $f$ is decreasing in the second variable.

Another example can be the following: $f(i,\eta_e)$ represents the time needed for clients of company $C_i$ to cross the edge $e$: the larger the number of  companies that own this edge is, and the larger the number of users to take this portion will be, and the longer the time needed to cross the edge will be. In this case, the function $f$ does not depend on $i$ and is increasing in the second variable.

The next result gives the existence of an asymptotic shape in this context. In other words,the travals times for a client of $C_i$ will give raise to a deterministic asymptotic shape that does not on the random environment (neither his company road network, nor the other networks), but only on the densities $(p_i)_{1 _le i \le n}$ of the distinct networks and on the function $f$. The same result is of course available for the cost point of view.

\begin{theorem}
Let $\Omega=(\Omega_S)^n$ and $\P=\P_{p_1}\times \P_{p_2} \dots \times \P_{p_n}$, with $p_i>p_c$ for each $i\in\{1,\dots,n\}$ and $f:\{1,\dots,n\}\times\{0,\dots,n\}\to\R_{+}$ with $f(i,n)>0$ as soon as $n\ne 0$.
We say that  two points $k$ and $l$ of $\Zd$ are $i$-connected
if there is a path between $k$ and $l$ such that $\omega_{e,i}=1$ for each bond $e$ which is used in the path. For $x\in\Zd$, we denote by $C_i(x)$ the $i$-connected component of $x$. 
Let us also define for $e=\{a,b\}\in\Ed$:
$$\eta_{e}=\sum_{i=1}^n \omega_{e,i}\1_{|C_i(a)|+|C_i(b)|=+\infty}$$  
As before, for $(\omega,\eta)\in\Omega$, and $(x,y)\in\Zd\times\Zd$, we define the
travel time
$d_i(x,y)(\omega,\eta)$ to be 
$$\inf_{\gamma}\sum_{e\in \gamma} f(i,\eta_e),$$
where the infimum is taken on the set of paths whose extremities are $x$
and $y$ and that are $i$-open in the configuration $\omega$.
For $t\ge 0$ and $i\in\{1,\dots,d\}$, we note
$$B_t^i=\{k\in\Zd; d_i(0,k)\le t\}.$$ 
For $i\in\{1,\dots,d\}$, let us note by $\Pcond_i$ the probability measure defined by  $$\Pcond_i(A)=\frac{\P(A\cap\{|C_i(0)|=+\infty\})}{\P(|C_i(0)|=+\infty)}.$$
Then, there is a convex compact set $A$ with non-empty interior such that
 $$\lim_{t\to +\infty}\frac{B_t}t =A \quad\Pcond_i\as$$
where the convergence holds for Hausdorff's topology. 
\end{theorem}    

\begin{proof}
Fix $i\in\{1,\dots,n\}$.
For $e=\{a,b\}\in\Ed$, let us define
$$\eta^{*}_{e}=f(i,1+\sum_{j\in\{1,\dots,n\}\backslash\{i\}} \omega_{e,j}\1_{|C_j(a)|+|C_j(b)|=+\infty})$$
and the
travel time
$d^{*}(x,y)(\omega,\eta)$ to be 
$$\inf_{\gamma}\sum_{e\in \gamma}\eta^{*}_e,$$
where the infimum is taken on the set of paths whose extremities are $x$
and $y$ and that are $i$-open in the configuration $\omega$.
For $t\ge 0$, we note
$$B_t^{*}=\{k\in\Zd; d^{*}(0,k)\le t\}.$$
It is clear that
$$\Pcond_i(\forall t\ge 0\quad B^i_t=B^{*}_t)=1.$$
Since the passage time $(\eta^{*}_e)_{e\in\Ed}$ are independent from the environment $(\omega_{k,i})_{k\in\Zd}$, it will be possible to use the preceding theorems.

The first step is to check that the law of  $(\eta^{*}_e)_{e\in\Ed}$ is
shift-invariant and ergodic under the translations of $\Zd$.
Obviously $$\miniop{}{\otimes}{j\in\{1,\dots,n\}\backslash\{i\}}\P_{p_j}=\big(\miniop{}{\otimes}{j\in\{1,\dots,n\}\backslash\{i\}}(1-p_j)\delta_0+p_j\delta_1\big)^{\otimes \Ed}$$ is shift-invariant and ergodic under the translations of $\Zd$. For simplicity, we note
$L=\big(\miniop{}{\otimes}{j\in\{1,\dots,n\}\backslash\{i\}}\delta_0+p_j\delta_1\big)^{\otimes \Ed}$.
We note 

$$F(\omega)=\big(f(i,1+\sum_{j\in\{1,\dots,n\}\backslash\{i\}} \omega_{e,j}\1_{|C_j(a(e))|+|C_j(b(e))|=+\infty})\big)_{e\in\Ed},$$
with the notation $e=\{a(e),b(e)\}$.

The law of  $(\eta^{*}_e)_{e\in\Ed}$ is precisely
the image of $L$ by the transformation $F$.
Let us denote by $S_{\nu}$ this law. 
It is easy to see that $F\circ \theta_u=\theta_u\circ F$ holds for each $u\in\Zd$, so that we have the commuting diagram:

\[
\begin{CD}
(\{0,1\}^{n-1})^{\Ed} @>\theta_u >> (\{0,1\}^{n-1})^{\Ed}\\
@VVFV @VVFV\\
(\R_{+})^{\Ed} @>\theta_u >> (\R_{+})^{\Ed}\\
\end{CD}
\]

For any event $A\in\bor[(\R_{+})^{\Ed}]$, we have
\begin{eqnarray*}S_{\nu}(\theta_u^{-1}(A)) & = & L(F^{-1}(\theta_u^{-1}(A))\\
 & = &  L( (\theta_u\circ F)^{-1}(A))\\
& = & L( (F\circ \theta_u)^{-1}(A))\\
& = & L(\theta_{u}^{-1}(F^{-1}(A))\\
& = & L(F^{-1}(A))=S_{\nu}(A).
\end{eqnarray*}
Then, $S_{\nu}$ is invariant under the translations.

It is now clear that the dynamical system $((\R_{+})^{\Ed},\bor[(\R_{+})^{\Ed}],S_{\nu},\theta_u)$
is a factor of the dynamical system $((\{0,1\}^{n-1})^{\Ed},\bor[(\{0,1\}^{n-1})^{\Ed}],L,\theta_u)$.
Since the dynamical system $((\{0,1\}^{n-1})^{\Ed},\bor[(\{0,1\}^{n-1})^{\Ed}],L,\theta_u)$ is ergodic,
it follows that the dynamical system $((\R_{+})^{\Ed},\bor[(\R_{+})^{\Ed}],S_{\nu},\theta_u)$ is ergodic too.

Now, since $\min f(i,1+.)>0$ and $\max f(i,1+.)<+\infty$, the conclusion follows from the remark that we  have done before the proof of theorem~\ref{asymptotic-shape}.
\end{proof}

\def\refname{References}
\bibliographystyle{alpha}
\bibliography{fpppc-final}

\end{document}